\documentclass[twocolumn]{IEEEtran}
\usepackage{amsmath,amssymb,amsthm,epsfig,color,empheq,graphicx,graphics,balance}
\usepackage{enumerate,url,wasysym,epstopdf,enumitem,array}
\usepackage{xcolor}
\usepackage{amsfonts}
\usepackage{adjustbox}
\usepackage{multirow}
\usepackage{cite}
\usepackage{accents}

\usepackage{algorithm2e} 
\usepackage{float} 
\DeclareMathOperator{\diag}{dg}

\DeclareMathOperator*{\argmin}{arg\,min}
\newtheorem{proposition}{Proposition}

\newtheorem{definition}{Definition}
\newtheorem{remark}{Remark}

\newcommand \bzero{\mathbf{0}}
\newcommand \bone{\mathbf{1}}
\newcommand \ba{\mathbf{a}}
\newcommand \bb{\mathbf{b}}
\newcommand \bc{\mathbf{c}}

\newcommand \bef{\mathbf{f}} % exception since \bf is taken

\newcommand \bp{\mathbf{p}}
\newcommand \bq{\mathbf{q}}

\newcommand \bv{\mathbf{v}}
\newcommand \bw{\mathbf{w}}

\newcommand \bz{\mathbf{z}}

\newcommand \bR{\mathbf{R}}

\newcommand \bX{\mathbf{X}}

\newcommand \balpha{\boldsymbol{\alpha}}

\newcommand \bdelta{\boldsymbol{\delta}}

\newcommand \blambda{\boldsymbol{\lambda}}
\newcommand \bmu{\boldsymbol{\mu}}

\newcommand \bsigma{\boldsymbol{\sigma}}

\newcommand \mcB{\mathcal{B}}

\newcommand \mcG{\mathcal{G}}

\newcommand \mcN{\mathcal{N}}

\newcommand \tbv{\tilde{\mathbf{v}}}

\newcommand \tbx{\tilde{\mathbf{x}}}

\newcommand \tbz{\tilde{\mathbf{z}}}

\newcommand \hbq{\hat{\mathbf{q}}}

\newcommand \hbx{\hat{\mathbf{x}}}

\newcommand \hbz{\hat{\mathbf{z}}}

\newcommand \bbq{\bar{\mathbf{q}}}

\newcommand \bbv{\bar{\mathbf{v}}}

\begin{document}

\title{Optimal Design of Volt/VAR Control Rules\\
of Inverters using Deep Learning}

\author{
    Sarthak Gupta,~\IEEEmembership{Graduate Student Member,~IEEE,}
    Vassilis Kekatos,~\IEEEmembership{Senior Member,~IEEE,}\\
    and
    Spyros~Chatzivasileiadis,~\IEEEmembership{Senior Member,~IEEE}
	
\thanks{Manuscript received March 30, 2023; revised September 19, 2023, and January 30, 2024; accepted March 15, 2024. Date of publication DATE; date of current version DATE. Paper no.TSG-00457-2023.}

\thanks{This work was supported by the US National Science Foundation (NSF) under grants 2034137 and 2150596, and by the European Research Council (ERC) Starting Grant VeriPhIED, Grant Agreement No. 949899. S. Gupta is with C3.AI, Redwood City, CA, USA. V. Kekatos is with the Elmore Family School of Electrical and Computer Engineering, Purdue University, IN 47906, USA. S. Chatzivasileiadis is with the Wind and Energy Systems Department, Technical University of Denmark, Lyngby, Denmark.}

\thanks{Color versions of one or more of the figures is this paper are available online at {http://ieeexplore.ieee.org}.}
\thanks{Digital Object Identifier XXXXXX}
}	
	
\markboth{IEEE TRANSACTIONS ON SMART GRID (to appear)}{Gupta, Kekatos, and Chatzivasileiadis: Optimal Design of Volt/VAR Control Rules of Inverters using Deep Learning}

\maketitle

\begin{abstract}
Distribution grids are challenged by rapid voltage fluctuations induced by variable power injections from distributed energy resources (DERs). To regulate voltage, the IEEE Standard 1547 recommends each DER inject reactive power according to piecewise-affine Volt/VAR control rules. Although the standard suggests a default shape, the rule can be customized per bus. This task of optimal rule design (ORD) is challenging as Volt/VAR rules introduce nonlinear dynamics, and lurk trade-offs between stability and steady-state voltage profiles. ORD is formulated as a mixed-integer nonlinear program (MINLP), but scales unfavorably with the problem size. Towards a more efficient solution, we reformulate ORD as a deep learning problem. The idea is to design a DNN that emulates Volt/VAR dynamics. The DNN takes grid scenarios as inputs, rule parameters as weights, and outputs equilibrium voltages. Optimal rule parameters can be found by training the DNN so its output approaches unity for various scenarios. The DNN is only used to optimize rules and is never employed in the field. While dealing with ORD, we also review and expand on stability conditions and convergence rates for Volt/VAR dynamics on single- and multi-phase feeders. Tests showcase the merit of DNN-based ORD by benchmarking it against its MINLP counterpart.
\end{abstract}
	

% \vspace*{-1em}

\section{Introduction}\label{sec:introduction}
DERs such as solar photovoltaics, are being advocated as a means to battle climate change, shave peak demand, and improve reliability. Despite the obvious benefits, the operation of distribution grids is nowadays challenged by undesirable voltage excursions induced by power injections from DERs. Traditional voltage regulation apparatus (e.g., regulators and capacitors) is deemed ineffective as responding to rapid and frequent voltage fluctuations can significantly shorten the lifetime of such equipment. Fortunately, the inverters interfacing DERs to the feeder can assist in regulating voltage by injecting reactive power. To this end, the IEEE Std. 1547 provisions that each DER should act autonomously, read its local voltage, and compute its reactive injection based on a Volt/VAR curve. This work aims to design these Volt/VAR curves optimally and customize them per bus according to the anticipated grid loading conditions every few hours or so.

Ideally, the reactive power setpoints of DERs can be optimally selected by solving an optimal power flow (OPF) given the current grid conditions~\cite{FCL,YGG13,8889689}. Such centralized schemes, however, involve high communication and computational overhead, privacy issues, and can introduce details. On the other hand, \emph{local} voltage regulation schemes for DERs entail calculating control setpoints solely based on locally collected data, such as load, solar generation, and voltage measurements at the grid interface. The Std. 1547 prescribes a local control scheme whereby DER setpoints are produced by \emph{control rules} taking the form of piecewise linear functions of local measurements~\cite{IEEE1547}. Local rules, though, are known to produce sub-optimal setpoints~\cite{Guido16,8667359}. Nevertheless, autonomy and simplicity are lucrative features of local schemes for real-time DER control. Focusing on Volt/VAR control, this work delves into the study and optimal design of local rules. Since voltages are affected by reactive setpoints, Volt/VAR rules give rise to closed-loop dynamics, which can become unstable under control rules with steep slopes~\cite{FCL13,7747940,9091863}. While the aforementioned works study the convergence and stability of Volt/VAR control rules, they do not address how to design such rules, i.e., how to select their exact shape, in the first place. Prior efforts on designing DER rules either resort to heuristics~\cite{6601722,6727491,8365842}, deal with non-dynamic Watt/VAR rules~\cite{9233917}; or restrict themselves to affine Volt/VAR rules~\cite{8003321}; or ignore the deadband~\cite{9781808}. A recent work designs optimal IEEE 1547-type rules, using mixed integer programs, which are then solved using relaxation heuristics~\cite{9609090}. However, discussions on the stability and convergence of the control rules are omitted. 

Similarly, there is a recent line of works that integrate Volt/VAR and Watt/VAR rules into an OPF formulation, and co-optimize rules and inverter setpoints for a single grid loading scenario~\cite{Paudyal1,Paudyal2,PaudyalMISOCP,Paudyal-LinDist3Flow}. Optimizing rules at the same timescale with setpoints may defy the intention of the IEEE Standard 1547 to have DERs either operating autonomously by running localized rules, or following centrally computed OPF setpoints. Regardless, references~\cite{Paudyal1,Paudyal2} co-optimize rules and setpoints via a nonlinear non-convex program, which may lack global optimality guarantees. In pursuit of global optimality, subsequent reference~\cite{PaudyalMISOCP} devised a mixed-integer second-order cone program (MISOCP), which could scale unfavorably with the network size (no running times were reported) and is limited to single-phase radial grids. Extending it to multiphase feeders under the exact AC grid model would call for computationally costly mixed-integer semidefinite programs (MISDP). As a remedy, reference~\cite{Paudyal-LinDist3Flow} adopts a linearized grid model to co-optimize rules and setpoints under the multiphase setting. Piecewise-affine rules are captured via binary variables selecting segments of voltages. Products of binary-continuous (voltage) variables are handled through the standard big-M trick to arrive at a mixed-integer linear program (MILP). However, the aforesaid works ignore linear inequality constraints on rule parameters required by IEEE Std. 1547 and do not ensure the rules are stable. Dealing with these two requirements requires keeping the native rule parameters (deadband, slope, saturation, reference voltage) as the optimization variables. Unfortunately, such parameterization of the problem introduces products between continuous variables, whose big-M reformulations are known not to be exact. This, in turn, gives rise to mixed-integer \emph{nonlinear} programs (MINLP). Reference~\cite{9656605} develops another MINLP formulation for ORD, which nonetheless, ignores voltage deadband, fixes kVAR saturation at the kVAR capacity of the inverter, does not enforce the constraints required by the IEEE Std. 1547, and is limited to single-phase grids.

Unlike above, our recent work in~\cite{MGCK23} formulates a bilevel optimization to design the slopes, deadband, saturation, and reference voltages for the Volt/VAR rules as the IEEE 1547.8 Standard prescribes. The bilevel optimization considers multiple grid scenarios to capture uncertainty. Upon leveraging the properties of the system at equilibrium, it finds stationary points using projected gradient descent iterates. 

Our present work extends~\cite{MGCK23}, and improves upon the previously cited literature in four directions:

\emph{c1)} Most of the existing works focus on simplified single-phase distribution grid models. We extend the analysis of~\cite{9091863}, and provide conditions to ensure the stability of the IEEE 1547 Volt/VAR rules for the more practically relevant setup of \emph{multiphase feeders};

\emph{c2)} We cast the problem of finding optimal Volt/VAR rules as a DNN training task. The training process involves stochastic projected gradient updates (SPGD) that leverage efficient, off-the-shelf Python libraries;

\emph{c3)} We genuinely design the DNN to emulate Volt/VAR dynamics: It accepts grid conditions as input, the parameters of the Volt/VAR rule as weights, and computes approximate equilibrium voltages at its output. Based on the convergence rate of Volt/VAR dynamics, we determine the minimum depth this DNN should have to approximate equilibrium voltages to the desired level of accuracy;

\emph{c4)} Leveraging the bilevel structure of ORD, we also reformulate ORD as a mixed-integer nonlinear program (MINLP). This MINLP-based approach does not scale favorably with the number of DERs and grid scenarios. Nonetheless, it serves as a benchmark for comparison to better assess the optimality and computational speed of our DNN-based ORD approach. 

We next expound upon how our work differs from prior works utilizing machine learning and/or reinforcement learning for smart inverter control. DNNs have been extensively employed before for optimal DER control under OPF formulations, with the objective of minimizing energy losses and energy costs; see e.g., \cite{Nanpeng19,ZhangISU,OPFandLearnTSG21,SGKCB2020}. Support vector machines and Gaussian processes have also been suggested for reactive power control using smart inverters~\cite{JKGD19,JSKGL22}. However, none of the above references aim at modeling the IEEE 1547-type piecewise linear, local, Volt/VAR rules. Furthermore, the existing problem formulations preclude the presence of closed-loop dynamics and are not nuanced by stability concerns as in the present work. In terms of using a DNN to model piecewise linear control rules, our work bears some similarities with references ~\cite{Baosen22,Baosen23,Low_ACC_22}. References~\cite{Baosen22} and \cite{Low_ACC_22} model piecewise linear control rule using a NN with a single hidden layer, and do not capture Volt/VAR dynamics over time using a neural network architecture. Furthermore, they focus on optimal control of transient dynamics. In contrast, the present work aims to design control rules that produce equilibrium voltages close to unity across many scenarios. This is achieved via efficient training of a recurrent neural network (RNN) whose training coincides with ORD. While the recent work~\cite{Baosen23} does leverage RNNs to design rules, it does so for controlling frequency transients.  Moreover, references~\cite{Baosen22,Baosen23,Low_ACC_22} do not discuss other topics covered in this work such as the IEEE 1547-type Volt/VAR rules, their convergence speed and depth of the resulting RNNs, as well as the implications of Volt/VAR control in multiphase feeders. Parallel work \cite{Cavraro22} designs a stable Volt/VAR control mechanism wherein each DER decides its reactive injections via a single-layer DNN driven by local data. Different from our approach, the suggested control rule does not comply with IEEE 1547 and training that DNN requires solving several OPF instances beforehand to generate labels. Finally, reference~\cite{Cavraro22} does not cover multiphase feeders.

\allowdisplaybreaks

%%%%%%%%%%%%%%%%%%%%%%%%%%%%%%%%%%%%%%%%%%%%%%%%%%%%%%%%%
\section{Feeder Modeling Preliminaries}\label{sec:model}
Consider a feeder rooted at the substation. Although the feeder can be single-phase or multiphase, it features a tree (radial) structure in terms of buses. For multiphase feeders, a bus may serve one to three phases; a valid pair of bus and phase will be referred to as a \emph{node}. For single-phase feeders, the terms \emph{bus} and \emph{node} will be used interchangeably. The substation is indexed by $0$ and is considered balanced; all remaining nodes are indexed by $n\in\mcN:=\{1,\dots,N\}$. All DERs are assumed to be single-phase and be able to provide reactive power control. For simplicity, each node is assumed to host a DER; we briefly discuss the minor modifications to deal with the more practical setting where not all nodes host DERs. Our numerical tests evaluate the latter setting. 

To study the effect of power injections on voltage magnitudes, we use an approximate linearized grid model. Let the active/reactive power injections and voltage magnitudes (henceforth simply voltages) at the non-substation nodes be collected into the $N$-length vectors $\bp$, $\bq$, and $\bv$, respectively. The linearized grid model relates these quantities as \cite{TJKT-SG21}
\begin{equation}\label{eq:ldf}
\bv\simeq \bR\bp+\bX\bq+v_0\bone    
\end{equation}
where $v_0$ is the substation voltage, and real-valued matrices $\bR$ and $\bX$ depend on line impedances and feeder topology.

If $\bp^g$ and $\bp^{\ell}$ denote the active power generated by DERs and that consumed by the loads accordingly, then $\bp = \bp^g - \bp^\ell$. Reactive power injections can be decomposed similarly as $\bq = \bq^g - \bq^\ell$. Supposing $\bp$ and $\bq^{\ell}$ are uncontrolled and vary with time, reactive power compensation entails adjusting $\bq^g$ to maintain $\bv$ around one per unit (pu). To isolate the effect of DER reactive injections on voltages, rearrange \eqref{eq:ldf} as
\begin{equation}\label{eq:ldf2}
\bv= \bX\bq^g+\tbv= \bX\bq+\tbv
\end{equation}
where the notation is slightly abused by denoting $\bq^g$ as $\bq$ for simplicity. The uncontrolled quantities are captured in vector $\tbv:=\bR(\bp^g - \bp^\ell)-\bX\bq^\ell+v_0\bone$, where $\tbv$ models voltages had it not been for reactive power compensation. Vector $\tbv$ will be henceforth termed the vector of \emph{grid conditions}. 

Given its importance in Volt/VAR control, let us summarize some properties of the sensitivity matrix $\bX$ appearing in \eqref{eq:ldf2}. For single-phase feeders, matrix $\bX$ is known to be symmetric, positive definite, and with positive entries; see e.g.,~\cite{BW2}, \cite{FCL13}. For multiphase feeders, however, matrix $\bX$ is non-symmetric and has positive as well as negative entries~\cite{VKZG16}. Nonetheless, under conditions typically met in practice~\cite{VKZG16}, matrix $\bX$ remains positive definite for multiphase feeders in the sense $\bz^\top\bX\bz>0$ for all $\bz\neq \bzero$. These nuances of $\bX$ call for relatively different treatments of Volt/VAR control between single- (Sections~\ref{sec:ni-single}--\ref{sec:1pDNN}) and multi-phase feeders (Section~\ref{sec:3pDNN}).

\begin{figure}[t]
	\centering
	\includegraphics[scale=0.63]{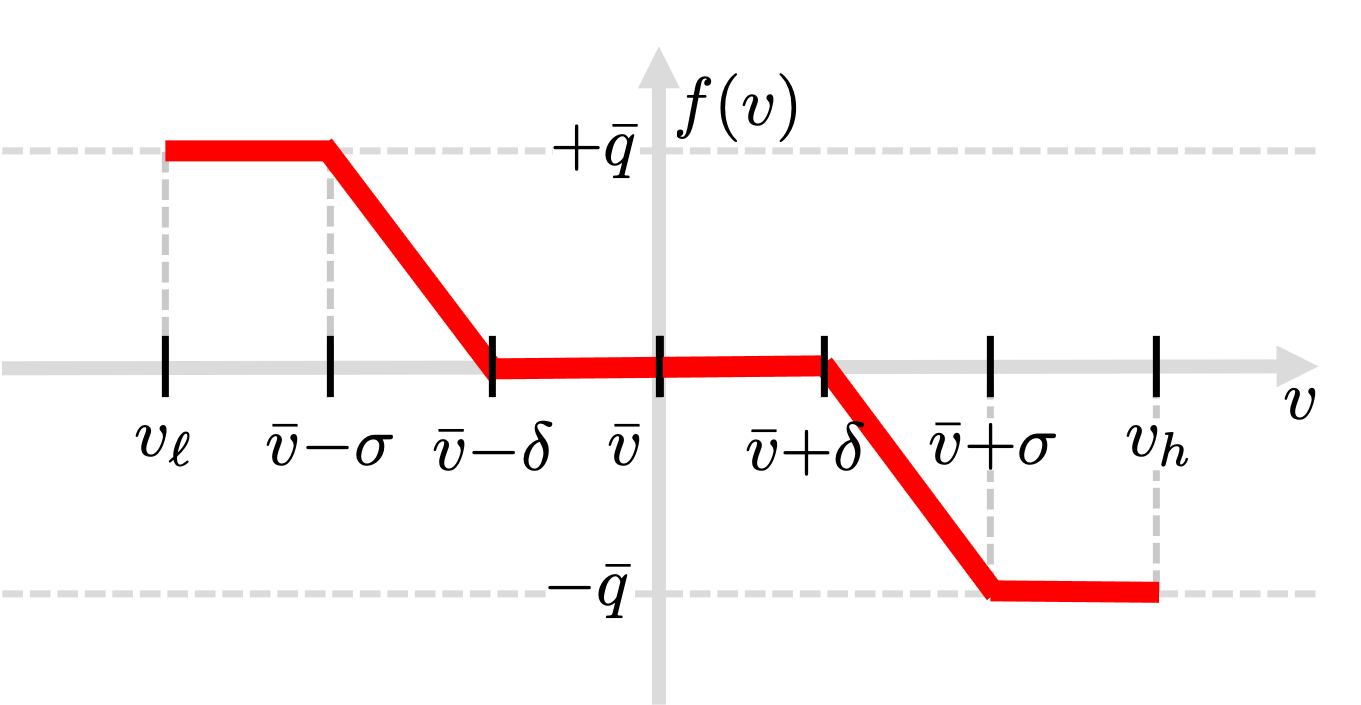}
% 	\vspace*{-1.5em}
	\caption{The piecewise linear Volt/VAR control rule $f(v)$ provisioned by the IEEE 1547 standard~\cite{IEEE1547}. The $x$-axis corresponds to the local voltage magnitude and the $y$-axis to the inverter setpoint for reactive power injection.}
% 	\vspace*{-1.5em}
	\label{fig:curve}
\end{figure}

% We consider inverter-interfaced DERs connected to a single phase, so their control rules remain the same for single-/multiphase feeders.

\section{Control Rules for Single-Phase Feeders}\label{sec:ni-single}
The IEEE 1547.8 standard provisions four modes of reactive power control~\cite{IEEE1547}: constant power, constant power factor, Watt/VAR, and Volt/VAR. We focus on the last one as being the most grid-adaptive. This mode enables the inverters to respond to local voltage deviations via a piecewise linear control curve $f(v)$, like the one depicted in Fig.~\ref{fig:curve}. The curve consists of a deadband of length $2\delta$ centered around $\bar{v}$; two affine regions; and two regions wherein reactive injections saturate at $\pm \bar{q}$. The IEEE standard constraints the curve parameters as follows (see Table~8 of~\cite{IEEE1547})
\begin{subequations}\label{eq:1547con2}
\begin{align}
0.95 &\leq \bar{v} \leq 1.05\label{eq:1547con2:v}\\
0 &\leq   \delta     \leq 0.03\label{eq:1547con2:delta}\\
\delta+0.02  &\leq   \sigma     \leq 0.18\label{eq:1547con2:sigma}\\
0&\leq \overline{q}\leq \hat{q}.\label{eq:1547con2:q}
\end{align}
\end{subequations}
Per \eqref{eq:1547con2:q}, the saturation value $\bar{q}$ can be equal to the reactive power capability $\hat{q}$ of the inverter, but also smaller than that. %Intuitively, each inverter is requested to inject (consume) reactive power when experiencing under-voltage (over-voltage); the deadband in the center prevents excessive switching between positive and negative injections; and saturation limits the slopes of linear segments. 

The rule of Fig.~\ref{fig:curve} is parameterized by $(\bar{v},\delta,\sigma,\bar{q})$, which can be customized per bus $n$ as $(\bar{v}_n,\delta_n,\sigma_n,\bar{q}_n)$. The rule can be alternatively parameterized by $(\bar{v}_n,\alpha_n,\delta_n,\bar{q}_n)$, where $\alpha_n$ is the negative slope of the affine segment and is defined as
\begin{equation}\label{eq:slope}
\alpha_n=\frac{\overline{q}_n}{\sigma_n-\delta_n}>0.
\end{equation}
Let vectors $(\bbv,\balpha,\bdelta,\bbq)$ collect $(\bar{v}_n,\alpha_n,\delta_n,\bar{q}_n)$ for all $n\in\mcN$; and stack such vectors together in vector $\bz:=(\bbv,\balpha,\bdelta,\bbq)$.

% can be smaller than or equal to  reactive power injection $\bar{q}$ remains within the operating capability $\hat{q}$ of the inverter. From Figure~\ref{fig:curve}, the underlying intention behind the standard may be inferred. If each of the inverters follows its Volt/VAR curve and injects (consumes) reactive power when facing under- (over-) voltage locally, the overall voltage profile of the distribution grid may improve. Moreover, the deadband in the center prevents excessive switching between positive and negative reactive power compensation. Beyond these intuitive insights, the next section presents existing and new results on the impact of Volt/VAR-interfaced DERs on single- and multi-phase grids. 

The interaction of Volt/VAR-controlled DERs with the grid results in the non-linear discrete-time dynamics
\begin{subequations}\label{eq:dynamics1}
\begin{align}
\bv^t &= \bX\bq^t + \tbv\label{eq:dyn_v}\\
\bq^{t+1} &= \bef_\bz(\bv^t)\label{eq:dyn_q}
\end{align}
\end{subequations}
where vector function $\bef_\bz(\bv^t)$ represents the action of Volt/VAR rules across all nodes and is parameterized by $\bz$. 

References~\cite{9091863,Pedram13} guarantee that Volt/VAR dynamics are stable if $\|\diag(\balpha)\bX\|_2<1$, where $\diag(\balpha)$ is a diagonal matrix having $\balpha$ on its diagonal. To be satisfied as a strict inequality, the condition can be strengthened as $\|\diag(\balpha)\bX\|_2\leq{1-\epsilon}$ for some $\epsilon\in(0,1)$. 

\begin{definition}\label{def:stable}
Volt/VAR rules satisfying $\|\diag(\balpha)\bX\|_2\leq{1-\epsilon}$ for $\epsilon\in(0,1)$ will be henceforth termed $\epsilon$-stable.
\end{definition}

To avoid the spectral norm condition $\|\diag(\balpha)\bX\|_2\leq 1-\epsilon$, we have previously proposed the polytopic restriction~\cite{MGCK23}:
\begin{subequations}\label{eq:stability2}
\begin{align}
\bX\balpha& \leq (1-\epsilon) \bone\label{eq:stability2:C1}\\
\alpha_n&\leq \frac{{1-\epsilon}}{\sum_{m\in\mcN}X_{nm}},\quad\forall n\in\mcN.\label{eq:stability2:C2}
\end{align}
\end{subequations}

% \begin{align}\label{eq:inner}
% \bq^*(\bz,\tbv)=\arg\min_{-\bbq\leq \bq\leq \bbq}~&F(\bq):= V(\bq)+C(\bq).
% \end{align}
% The two components of the objective $F(\bq)$ are defined as
% \begin{subequations}\label{eq:costs}
% \begin{align}
% V(\bq)&:=\label{eq:Vn}\\
% C(\bq)&:=.\label{eq:Cn}
% \end{align}
% \end{subequations}

If stable, the dynamics in \eqref{eq:dynamics1} enjoy an equilibrium under any grid condition $\tbv$~\cite{FCL13}. In fact, the inverter setpoints at equilibrium coincide with the unique minimizer of the convex program~\cite{FCL13}
\begin{equation}\label{eq:inner}
\bq^*(\bz,\tbv):=\argmin_{-\bbq\leq \bq\leq \bbq}
F(\bq)
\end{equation}
where the objective function is defined as
\[F(\bq):=
\underbrace{\tfrac{1}{2}\bq^\top \bX\bq+\bq^\top(\tbv-\bbv)}_{:=V(\bq)}+\underbrace{\sum_{n\in\mcN}\tfrac{1}{2\alpha_n}q_n^2 +\delta_n|q_n|}_{:=C(\bq)}.\]
Component $V(\bq)$ can be equivalently expressed as~\cite{FCL13}
\begin{align}\label{eq:Vn_2}
   V(\bq)=\frac{1}{2}(\bv-\bbv)^\top \bX^{-1}(\bv-\bbv) +\mathrm{constants.} 
\end{align}

Because $\bX\succ 0$, function $V(\bq)$ is an $\ell_2$-norm of $(\bv-\bbv)$. Hence, minimizing $V(\bq)$ aims at bringing voltages close to reference voltages. Nonetheless, problem \eqref{eq:inner} involves also $C(\bq)$ in its cost. Hence, to best regulate voltages, one would try setting $\balpha$ to infinity and $\bdelta$ to zero so $C(\bq)=0$ and the equilibrium setpoints minimize only $V(\bq)$. This course of action however would violate the stability condition of $\|\diag(\balpha)\bX\|_2\leq 1$.

% The existence of the equilibrium alone does not ensure that the dynamics in \eqref{eq:dynamics1} are stable. 

The next section develops methods for selecting the Volt/VAR rule parameters $\bz$ so that a voltage regulation objective is minimized for a set of grid scenarios. For single-phase feeders, Section~\ref{sec:1pDNN} reformulates (ORD) as the problem of training a neural network, while Section~\ref{sec:1pMINLP} tackles ORD as a mixed-integer nonlinear program. For multiphase feeders, solving ORD is dealt with in Section~\ref{sec:3pDNN}.

%%%%%%%%%%%%%%%%%%%%%%%%%%%%%%%%%%%%%%%%%%%%%%%%%%
\section{ORD for $1\phi$ Feeders via Deep Learning}\label{sec:1pDNN}
Because Volt/VAR rules are used so inverters can operate autonomously, it is reasonable to assume that rule parameters $\bz$ are updated infrequently, say every 2 hours. Then, rules $\bz$ should be optimized while considering the possibly diverse loading conditions the feeder may experience over those 2 hours. To account for such conditions, suppose we are given a set of $S$ load/solar scenarios $\{(\bp_s^g,\bp_s^\ell,\bq_s^\ell)\}_{s=1}^S$. Each scenario is related to grid condition vector [see~\eqref{eq:ldf}]
\[\tbv_s:=\bR(\bp^g_s-\bp^\ell_s)-\bX\bq^\ell_s.\]
Let $\bq^*(\bz,\tbv_s)$ or simply $\bq^*_s(\bz)$ denote the equilibrium setpoints reached by stable Volt/VAR rules parameterized by $\bz$ under grid conditions $\tbv_s$. Unfortunately, setpoints $\bq_s^*(\bz)$ cannot be expressed as in closed form. They can be computed by either iterating \eqref{eq:dynamics1}, or as the minimizer of \eqref{eq:inner}. The related equilibrium voltage is $\bv_s^*(\bz):=\bX\bq^*_s(\bz)+\tbv_s$ from \eqref{eq:ldf}. 

We pose the ORD task as a minimization problem over $\bz$:
 \begin{align}\label{eq:outer}
\min_{\bz}~&~\frac{1}{2S} \sum_{s=1}^S \|\bX\bq^*_s(\bz) + \tbv_s-\bone\|_2^2\\
\mathrm{s.to}~&~\eqref{eq:1547con2},\eqref{eq:stability2}\nonumber
\end{align}
to minimize the Euclidean distance of equilibrium voltages from unity, averaged across scenarios. Constraints \eqref{eq:1547con2} and \eqref{eq:stability2} ensure rules are stable and compliant with the IEEE 1547. 

It is worth iterating that the ORD is solved centrally by the utility operator, but only every two hours or so. Once decided, the optimal rule parameters $\bz$ are communicated to the DERs, which can operate autonomously for the next two hours. The frequency at which the operator chooses to re-optimize the rules or the number of scenarios $S$ used, do not alter the proposed methodology per se. Increasing $S$ increases the running time to some moderate extent as demonstrated in the numerical tests. Re-optimizing the rules more frequently (say every one instead of every two hours) would apparently yield better grid performance. Nonetheless, it would raise the communication and computational cost for the utility. If the operator is willing to communicate with the DERs every 15 minutes or less, it may be more meaningful to solve a (stochastic) OPF and communicate direct setpoints rather than rules to DERs.

An additional important observation here is that \eqref{eq:outer} aims at minimizing voltage deviations from unity averaged across buses and scenarios $S$. Interestingly, the conference offshoot of this work shows how the proposed ORD methodology can be extended to design rules that minimize ohmic losses subject to voltages lying within the desired range~\cite{NAPS2023}. To account for uncertainty, losses are averaged over scenarios, and voltage ranges are enforced as chance constraints.

One may wonder why we are not satisfied with the fact that any stable rule $\bz$ settles at the minimizer of \eqref{eq:inner}, which is seemingly a meaningful equilibrium. Such equilibrium may be insufficient due to three reasons: \emph{i)} The term $V(\bq)$ is a \emph{rotated} $\ell_2$-norm of $(\bv-\bbv)$, so that voltage deviations are weighted unequally across buses; \emph{ii)} If DERs are sited only on a subset $\mcG\subset \mcN$ of nodes, the cost $V(\bq)$ gets modified as $V_{\mcG}(\bq_{\mcG})=\frac{1}{2}(\bv_{\mcG}-\bbv_{\mcG})^\top\bX_{\mcG\mcG}^{-1}(\bv_{\mcG}-\bbv_{\mcG})$, where subscript $\mcG$ denotes the subvectors/submatrix obtained by keeping the rows/columns corresponding to $\mcG$; see~\cite{MGCK23}. Such cost may not be representative of $\|\bv-\bbv\|_2^2$; and \emph{iii)} As discussed earlier, stability limitations do not allow us to set $\balpha$ to infinity although it seems desirable from a voltage regulation standpoint. The aforementioned reasons motivate the need to optimally design $\bz$ so the induced equilibrium voltages $\bv_s^*(\bz)$ are better regulated.

% \[V_{\mcG}=\frac{1}{2}\bq_{\mcG}^\top \bX_{\mcG\mcG}\bq_{\mcG}+\bq_{\mcG}^\top(\tbv_{\mcG}-\bbv_\mcG)=(\bv_{\mcG}-\bbv_{\mcG})^\top\bX_{\mcG\mcG}^{-1}(\bv_{\mcG}-\bbv_{\mcG})\]

Albeit simply stated, problem \eqref{eq:outer} is computationally challenging as $\bq_s^*(\bz)$ is the solution of the \emph{inner} minimization problem \eqref{eq:inner}, which is parameterized by $\bz$. Thereby, the ORD task is a \emph{bilevel optimization} over $\bz$: The \emph{outer} problem \eqref{eq:outer} depends on $S$ inner problems of the form \eqref{eq:inner}, one per scenario.

Our first strategy towards tackling \eqref{eq:outer} is to replace the inner problem with a DNN that simulates the Volt/VAR dynamics. This DNN has $\bz$ as weights, accepts $\tbv_s$ as input, and outputs the equilibrium voltages $\bv^*_s(\bz)$. Let the DNN output be denoted by $\Phi(\tbv_s;\bz)$. The key idea is that if $\Phi(\tbv_s;\bz)$ are the equilibrium voltages for rule $\bz$ over scenario $s$, then problem \eqref{eq:outer} becomes the supervised training task:
\begin{align}\label{eq:outer_DNN}
\min_{\bz}~&~L(\bz):=\frac{1}{2S} \sum_{s=1}^S \|\Phi\left(\tbv_s;\bz\right)-\bone\|_2^2\\
\mathrm{s.to}~&~\eqref{eq:1547con2},\eqref{eq:stability2}.\nonumber
\end{align}
To draw a useful analogy, scenarios $\tbv_s$ are analogous to feature vectors in regression problems; equilibrium voltages $\bv^*_s=\Phi(\tbv_s;\bz)$ are the predictions for feature vectors; and $\bone$ is the (constant) target label for the prediction. Formulating~\eqref{eq:outer} as~\eqref{eq:outer_DNN} allows us to leverage efficient DNN libraries for optimizing $\bz$. With this motivation in mind, we next design the DNN, and then describe the steps to train it.

%%%%%%%%%%%
\subsection{Designing a Digital Twin of Volt/VAR Dynamics}\label{subsec:design}
\begin{figure}[t]
	\centering
	\includegraphics[scale=0.65]{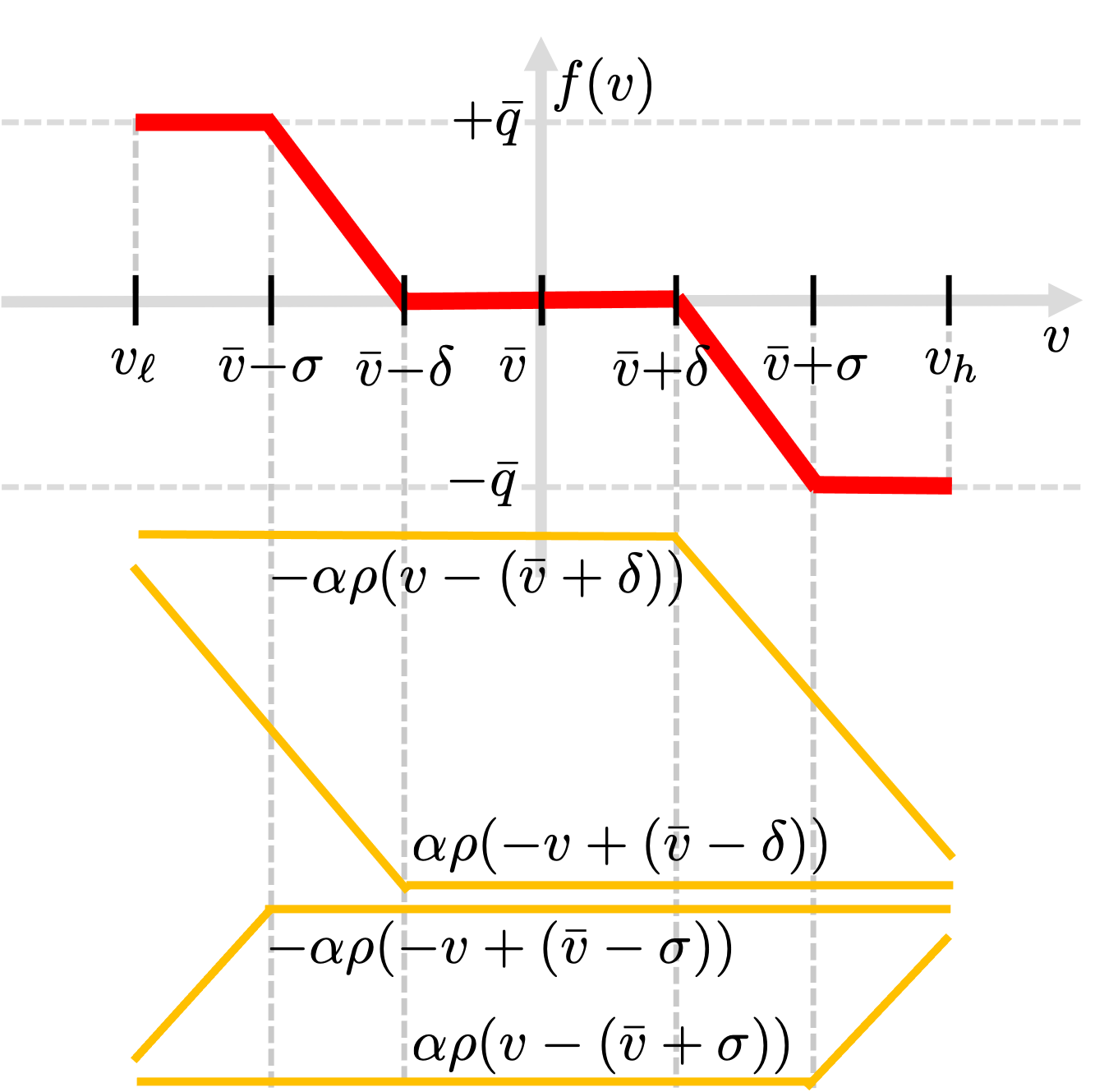}
% 	\vspace*{-1.5em}
	\caption{Volt/VAR rule $f(v)$ expressed as a sum of ReLUs.}
% 	\vspace*{-1.5em}
	\label{fig:relus}
\end{figure}

\begin{figure}[t]
\centering
	\includegraphics[scale=0.5]{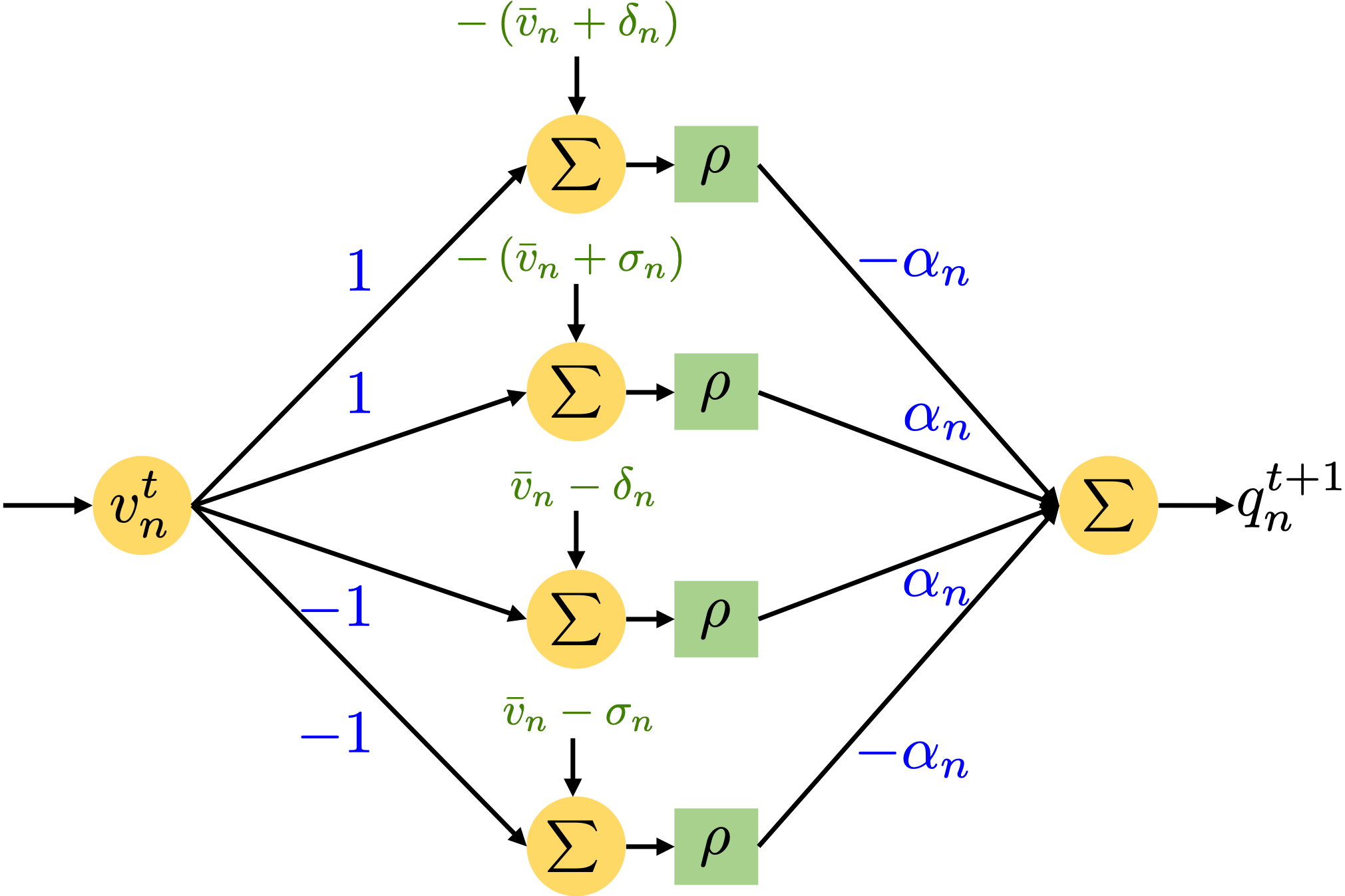}
% 	\vspace*{-1.5em}
	\caption{Volt/VAR rule $f(v)$ model using a NN with 1 hidden layer.}
% 	\vspace*{-1.5em}
	\label{fig:single}
\end{figure}

The Volt/VAR curve of Fig.~\ref{fig:curve} can be interpreted as a superposition of four piecewise-linear functions, each with a single breakpoint, as shown in Fig.~\ref{fig:relus}. These functions can be thought of as the outputs of rectified linear units (ReLU) $\rho(x)$, which return $x$ for $x>0$; and $0$ otherwise. To get the different breakpoints and slopes as in Figure~\ref{fig:relus}, the ReLU units need the appropriate inputs and scaling. The required mathematical operations can be implemented through the DNN of Fig.~\ref{fig:single}, which takes $v^t_n$ as input and computes the setpoint $q^{t+1}_n$ at its output. The input and output layers have one neuron each. The hidden layer consists of four neurons. The weights of the hidden layer are fixed to $[1,1,-1,-1]^{\top}$, but its bias vector is trainable and given by $[-(\bar{v}+\delta),-(\bar{v}+\sigma),\bar{v}-\delta,\bar{v}-\sigma]^{\top}$. Each of the four neurons in the hidden layer is equipped with a ReLU unit. The output layer has a trainable weight vector $[-\alpha,\alpha,\alpha,-\alpha]^{\top}$ and the bias is fixed at $0$. 

\begin{figure*}
	\centering
	\includegraphics[scale=0.44]{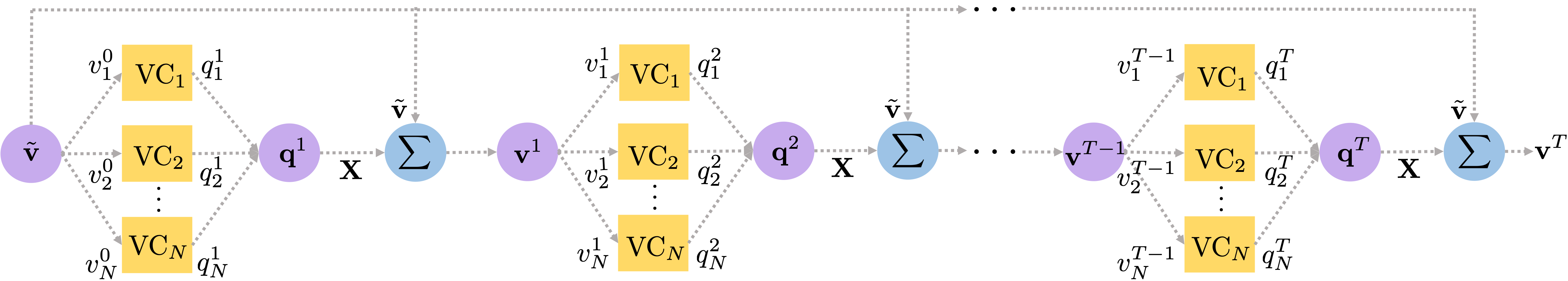}
% 	\vspace*{-1.5em}
	\caption{DNN-based digital twin for the Volt/VAR dynamics of~\eqref{eq:dynamics1}. The DNN is structured so that $T$ time steps are arranged horizontally. The modules $\text{VC}_n$'s implementing the Volt/VAR curves for each one of the $N$ inverters are stacked vertically. Skip connections propagate the input vector (grid scenario) $\tbv$ to each time instant to implement $\bv^{t+1}=\bX\bq^{t+1}+\tbv$.}
% 	\vspace*{-1.5em}
	\label{fig:DNN}
\end{figure*}

Heed that the NN of Fig.~\ref{fig:single} implements the Volt/VAR curve for a single inverter and a single time step as $q_n^{t+1}=f_n(v_n^t)$. To simulate the entire Volt/VAR network dynamics of \eqref{eq:dynamics1}, we will treat the NN of Fig.~\ref{fig:single} as a building block and replicate it across inverters and time. Let $\text{VC}_n$ represent the NN module running a single time step for inverter $n$. This module is parameterized by $(\bar{v}_n,\alpha_n,\sigma_n,\delta_n)$. With a slight abuse of terminology, let the collection of $\text{VC}_n$'s for all inverters be labeled as a single \emph{layer}. These modules are stacked vertically, as shown in Fig.~\ref{fig:DNN}. Each one of these layers implements \eqref{eq:dyn_q} by receiving $\bv^t$ as input and producing setpoints $\bq^{t+1}$ as output at time $t$. The setpoints $\bq^{t+1}$ in turn produce voltages $\bv^{t+1}=\bX\bq^{t+1}+\tbv$ per the grid model~\eqref{eq:dyn_v}. To simulate the dynamics over time, the new voltages $\bv^{t+1}$ are passed to the next layer corresponding to time $t+1$. The process is repeated for $T$ steps. Structurally, these interactions result in a larger DNN with $T$ repeating layers, one layer per iteration of the dynamics in \eqref{eq:dynamics1}, as shown in Figure \ref{fig:DNN}. The simulation of dynamics over $T$ iterations is equivalent to performing a forward pass through the larger DNN with $\tbv$ as the input. To implement \eqref{eq:dyn_v}, the input $\tbv$ (grid scenario) is also propagated to the inner layers via so-called \emph{skip connections}. 

\begin{figure}[t]
	\centering
	\includegraphics[scale=0.44]{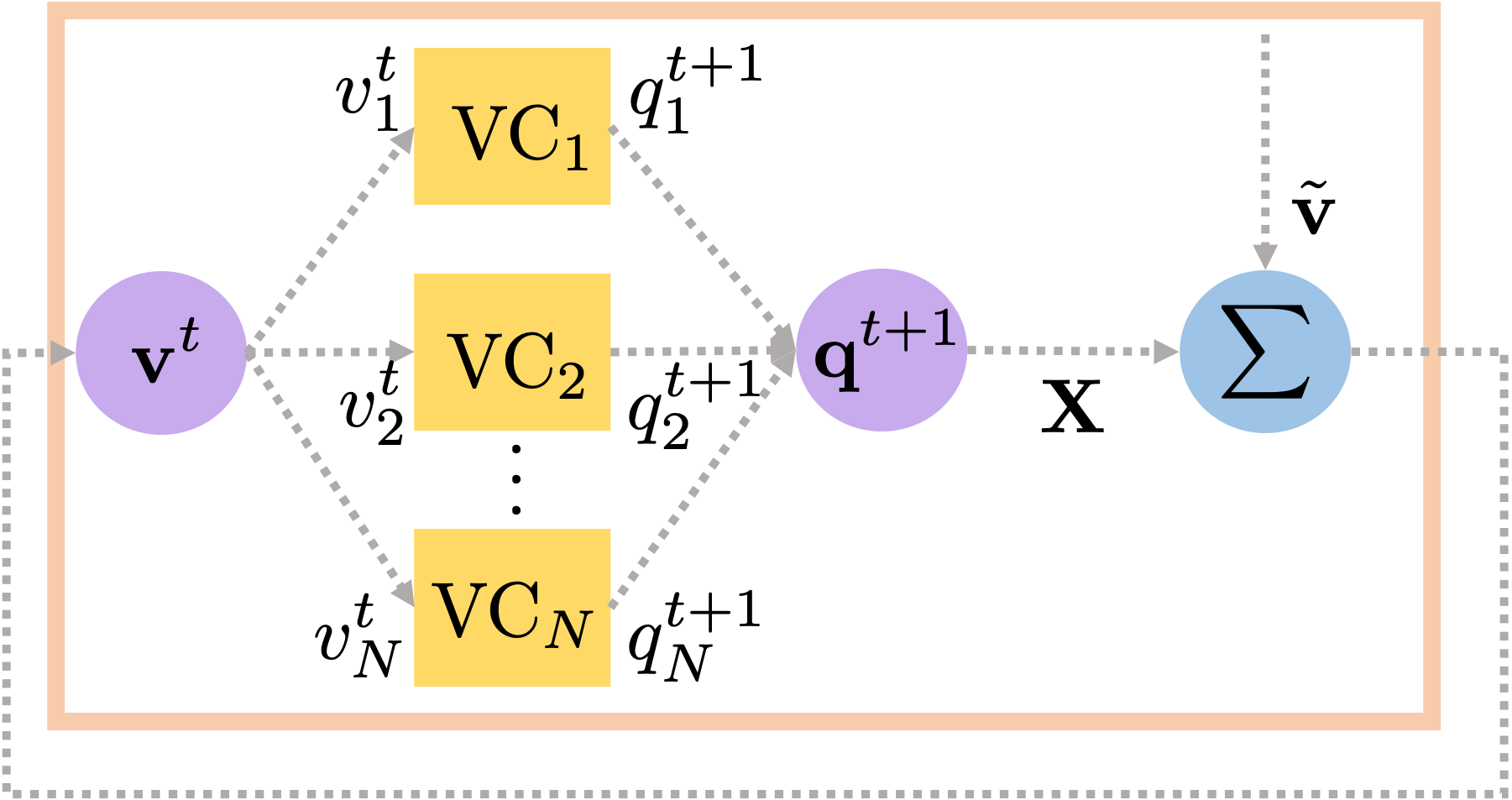}
% 	\vspace*{-1.5em}
	\caption{Recurrent representation (RNN) of the digital twin of Fig.~\ref{fig:DNN}.}
% 	\vspace*{-1.5em}
	\label{fig:RNN}
\end{figure}

It is worth stressing that each module $\text{VC}_n$ is replicated horizontally across the $T$ times. This implies significant \emph{weight sharing} across the $T$ layers. Therefore, the number of trainable parameters $\hbz:=(\bbv,\balpha,\bdelta,\bsigma)$ remains fixed at $4N$, irrespective of the DNN depth $T$. This weight-sharing aspect results in computational and memory-related efficiencies for DNN storing, prediction, and training, and has been instrumental in the success of architectures such as recurrent (RNN), convolutional (CNN), or graph (GNN) neural networks. In fact, it is possible to obtain a recurrent \emph{`rolled'} representation of the larger DNN of Fig.~\ref{fig:DNN}, as shown in Fig.~\ref{fig:RNN}, allowing one to utilize RNN-specific functionalities in DNN libraries.

In a nutshell, the DNN of Fig.~\ref{fig:DNN} simulates Volt/VAR dynamics across $T$ times. In other words, once fed with a grid condition vector $\tbv_s$, its output will approximate the equilibrium voltages $\bv_s^*$ reached by Volt/VAR dynamics under rule $\bz$. As a result, optimizing over $\bz$ by training the DNN so that equilibrium voltages $\{\bv_s^*(\bz)\}_{s=1}^S$ come close to one pu, serves the purposes of ORD. Surrogating Volt/VAR dynamics by the DNN is effective only if the DNN depth $T$ is sufficiently large. \emph{How deep should the DNN be so that its output $\Phi\left(\tbv_s;\hbz\right)$ is close to $\bv_s^*$?} Because a DNN of depth $T$ simulates exactly the Volt/VAR dynamics up to time $T$, the answer for selecting $T$ is apparently the settling time of the Volt/VAR dynamics as detailed next and shown in the appendix.

\begin{proposition}\label{pro:T}
Suppose $\epsilon$-stable Volt/VAR rules are described by $\bz$. The depth $T$ of the DNN in Fig.~\ref{fig:DNN} required to ensure $\|\Phi\left(\tbv;\bz\right) -\bv^*(\bz)\|_2\leq \epsilon_1$ for all grid conditions $\tbv$ is
\begin{equation*}%\label{eq:pro1}
T\geq \frac{ \log\frac{2\|\bX\|_2\|\hbq\|_2}{\epsilon_1}}{\log(1-\epsilon)^{-1}}.
\end{equation*}
\end{proposition}

The result implies that the minimum depth $T$ grows logarithmically with the desired accuracy $\epsilon_1$ and the stability margin $\epsilon$. Plugging in the typical values $\epsilon_1=10^{-4}$, $\|\bX\|_2=4.63\cdot10^{-1}$ for IEEE 37-bus feeder, $\|\hbq\|_2=0.1$, and $\epsilon=0.3$, the bound yields a comfortably small number of $T\geq 20$ layers. For $\epsilon_1=10^{-6}$, the number of layers $T$ increases to $32$, demonstrating the scalability of the approach. 

To recapitulate, the Volt/VAR rule for each DER is described by four parameters $(\bar{v}_n,\sigma_n,\delta_n,\alpha_n)$. These parameters appear as weights of the single-layer neural network shown in Fig.~\ref{fig:single}. This building block is indicated as a yellow block and termed $\text{VC}_n$ in the RNN of Figure~\ref{fig:DNN}. Each $\text{VC}_n$ is repeated at each layer $t$ of the RNN. Although the RNN may have $T$ layers, there is significant weight sharing and only $4N$ parameters to be trained, 4 per inverter.

%%%%%%%%%%%%%
\subsection{DNN Training}\label{subsec:training}
With rule parameters $\hbz$ embedded as DNN weights and biases, the optimal Volt/VAR curves are obtained by training $\Phi\left(\tbv_s;\hbz\right)$. Conventional DNN training uses stochastic gradient descent (SGD) to update the DNN parameters and eventually minimize the \emph{loss function} in~\eqref{eq:outer_DNN}. However, parameters $\hbz$ should satisfy constraints~\eqref{eq:1547con2} and \eqref{eq:stability2}. Plain SGD may fail to return a feasible $\bz$. This can be circumvented by using projected stochastic gradient (PSGD) updates. PSGD updates first compute an intermediate quantity $\hbx^{i+1}$ via gradient descent%~\cite{DL_book}
\begin{align}\label{eq:spgd1}
\hbx^{i+1}&=\hbz^{i}-\frac{\mu}{2B}{\nabla_{\hbz^i}}\left(\sum_{s\in \mcB_i} \|\Phi(\tbv_s;\hbz)-\bone\|_2^2\right)
\end{align}
where $\mu>0$ is the step size; set $\mcB_i$ is a batch of $B$ scenarios (a subset of the original $S$ scenarios); and $\nabla_{\hbz^i}(\cdot)$ is the gradient of the loss function with respect to $\hbz$ evaluated at $\hbz=\hbz^i$. The gradient term in \eqref{eq:spgd1} is calculated efficiently thanks to \emph{gradient back-propagation}.

The second step for PSGD updates entails projecting $\hbx^{i+1}$ into the feasible space defined by \eqref{eq:1547con2} and \eqref{eq:stability2}. To this end, we first transform $\hbx^{i+1}$ from parameter space $(\bbv,\balpha,\bdelta,\bsigma)$ to space $(\bbv,\bc,\bdelta,\bsigma)$, where vector $\bc$ has entries $c_n:=1/\alpha_n$. Variable $\hbx^{i+1}$ transformed in the new space is called $\tbx^{i+1}$. The transformation is a one-to-one mapping between the two spaces and is used so that the feasibility set induced by \eqref{eq:1547con2} and \eqref{eq:stability2} is convex, and so it is easy to project onto it. We proposed this transformation in \cite{MGCK23}. We review it here for completeness. Using \eqref{eq:slope}, constraint \eqref{eq:1547con2:q} is expressed as
\begin{align}\label{eq:ieee1547:qtilde}
    0\leq\bsigma-\bdelta\leq\bc\odot\hbq
\end{align}
where $\odot$ means element-wise multiplication. Constraints~\eqref{eq:stability2} can be expressed in terms of $\bc$ instead of $\balpha$ as~\cite{MGCK23}
\begin{subequations}\label{eq:stability3}
\begin{align}
\bc&\geq \frac{1}{{1-\epsilon}}\bX\bone\label{eq:stability3:a}\\
\bX\ba&\leq (1-\epsilon)\cdot \bone\label{eq:stability3:b}\\
\ba\odot \bc&\geq \bone,\quad \forall~n\in\mcN\label{eq:conic}
\end{align}
\end{subequations}
where $\ba$ is an auxiliary variable. Constraint~\eqref{eq:conic} can be rewritten as a second-order cone. In~\cite{MGCK23}, we show how \eqref{eq:stability3} is equivalent to \eqref{eq:stability2}. The quantity $\tbx^{i+1}$ can now be projected onto the feasible space via the convex minimization
\begin{align}\label{eq:project}
\tbz^{i+1} =\arg\min_{\bz}~&~\|\tbx^{i+1}-\bz\|_2^2\\
\mathrm{s.to}~&~\eqref{eq:1547con2:v}-\eqref{eq:1547con2:sigma},\eqref{eq:ieee1547:qtilde},\eqref{eq:stability3}.\nonumber
\end{align}
The PSGD update is completed by transforming $\tbz^{i+1}$ from space $(\bbv,\bc,\bdelta,\bsigma)$ back to space $(\bbv,\balpha,\bdelta,\bsigma)$ to get $\hbz^{i+1}$.

The proposed DNN training can be implemented in Python using DNN libraries such as PyTorch. Step~\eqref{eq:spgd1} is the standard SGD update pertaining to the loss function of~\eqref{eq:outer_DNN} over the batch of training labels $\{\tbv_s,1\}_{\mcB_i}$. As with standard DNN training, adaptive moment-based algorithms such as Adam can enable fast convergence and avoid saddle points. The DNN weights and biases are transformed between the parameter spaces and then passed to a convex optimization module to implement the projection step of~\eqref{eq:project}. In the last step, DNN weights and biases are updated with the new projected parameters, upon transformation to the original space. The steps are repeated for several epochs.
%~\cite{Adam}

%%%%%%%%%%%%%%%%%%%%%%%%%%%%%%%%
\section{ORD for $1\phi$ Feeders as an MINLP}\label{sec:1pMINLP}
A second approach towards solving the bilevel program in~\eqref{eq:outer} is to replace each inner problem with its first-order optimality conditions and append these conditions as constraints to the outer problem. To capture complementary slackness, we will introduce binary variables and use the so-termed big-M trick to eventually express the outer problem as a mixed-integer nonlinear program (MINLP). The process is delineated next. Although this MINLP approach does not scale gracefully with the number of DERs and/or scenarios, it serves as a benchmark for the DNN-based ORD.

We first transform \eqref{eq:inner} to a differentiable form as
\begin{subequations}\label{eq:inner2}
\begin{align}
\min_{\bq,\bw}~&~\frac{1}{2}\bq^\top \left(\bX+\text{dg}(\bc)\right)\bq+\bq^\top(\tbv_s-\bbv)+\bdelta^{\top}\bw\\
\mathrm{s.to~}&-\bw \leq \bq \leq\bw: \quad\quad\quad (\underline{\blambda},\overline{\blambda}) \label{eq:inner2:s}\\
~&~-\bbq \leq \bq \leq\bbq: \quad\quad\quad (\underline{\bmu},\overline{\bmu})\label{eq:inner2:qbar}
\end{align}
\end{subequations}
where vector $\bc$ has entries $c_n:=1/\alpha_n$, and variable $\bw$ has been introduced to deal with the non-differentiable terms $|q_n|$ in \eqref{eq:inner}. Slightly abusing notation, denote the optimal primal/dual variables of \eqref{eq:inner2} by $(\bq,\bw;\underline{\blambda},\overline{\blambda},\underline{\bmu},\overline{\bmu})$. Although the variables vary per scenario, we suppress subscript $s$ for simplicity. These variables satisfy the optimality conditions\begin{subequations}\label{eq:kkt}
\begin{align}
    \left(\bX+\diag(\bc)\right)\bq+\tbv_s-\bbv-\underline{\blambda}+\overline{\blambda}-\underline{\bmu}+\overline{\bmu}&=\bzero\label{eq:kkt:loq}\\
    \bdelta-\underline{\blambda}-\overline{\blambda}&=\bzero\label{eq:kkt:los}\\
    -\bw\leq \bq &\leq\bw\label{eq:kkt:pfs}\\
    -\bbq\leq \bq &\leq\bbq\label{eq:kkt:pfq}\\
    \underline{\blambda},\overline{\blambda},\underline{\bmu},\overline{\bmu}&\geq\bzero\label{eq:kkt:df}\\
    \overline{\blambda}\odot\left(\bq-\bw\right)&=\bzero\label{eq:kkt:cs1}\\
    \underline{\blambda}\odot\left(-\bq-\bw\right)&=\bzero\label{eq:kkt:cs2}\\
    \overline{\bmu}\odot\left(\bq-\bbq\right)&=\bzero\label{eq:kkt:cs3}\\
    \underline{\bmu}\odot\left(-\bq-\bbq\right)&=\bzero.\label{eq:kkt:cs4}
\end{align}
\end{subequations}
Equalities \eqref{eq:kkt:loq}--\eqref{eq:kkt:los} follow from Lagrangian optimality; and inequalities \eqref{eq:kkt:pfs}--\eqref{eq:kkt:df} from primal/dual feasibility. Inequalities \eqref{eq:kkt:cs1}--\eqref{eq:kkt:cs4} are complementary slackness conditions.

The bilevel problem in \eqref{eq:outer} can be now reduced to a single-level formulation upon appending conditions \eqref{eq:kkt} as constraints to \eqref{eq:outer} per scenario $s$. Such constraints ensure that $\bq_s^*(\bz)$ is indeed the minimizer of \eqref{eq:inner}. Nonetheless, constraint \eqref{eq:kkt:loq} and the complementary conditions introduce bilinear terms. Bilinearity can be partially addressed by handling complementary slackness conditions through the big-M trick. For example, condition \eqref{eq:kkt:cs1} can be expressed as
\begin{subequations}\label{eq:bigM}
\begin{align}
    \bzero\leq\overline{\blambda}&\leq M_1 \bb\\
    \bzero\leq\bq-\bw&\leq M_2 (\bone-\bb)
\end{align}
\end{subequations}
where $\bb$ is an $N$-dimensional binary variable, and $(M_1,M_2)$ are large positive constants. The latter can be selected as $M_2=2\bbq$, while the former can be set to a numerically estimated upper bound of the corresponding dual variables $\underline{\blambda}$. Although complementarity constraints can be reformulated to linear ones, that is not true for the bilinear term $\diag(\bc)\bq$ between two continuous variables appearing in \eqref{eq:kkt:loq}. This term gives rise to a mixed-integer nonlinear program rather than a mixed-integer linear program.

Since \eqref{eq:kkt}--\eqref{eq:bigM} contain $\bc$ and $\bbq$, the constraints \eqref{eq:1547con2} and \eqref{eq:stability2} need to be rewritten in terms of $\bc$ and $\bbq$ as well. To this end, we chose the parameterization $\tbz:=(\bbv,\bc,\bdelta,\bbq)$. In this new parameterization, constraint~\eqref{eq:1547con2:sigma} is replaced by
\begin{align}
0.02\cdot \bone  &\leq \bc\odot\bbq    \leq 0.18\cdot \bone-\bdelta\label{eq:1547con3:cq}
\end{align}
which introduces bilinear terms too. Stability constraints~\eqref{eq:stability2} have already been transformed from $\balpha$ to $\bc$ in \eqref{eq:stability3}. 

Putting everything together, the bilevel ORD problem of \eqref{eq:outer} can be solved as the MINLP:
\begin{subequations}\label{eq:MINLP}
 \begin{align}
\tbz^*=\arg\min_{\tbz}~&~\frac{1}{2S} \sum_{s=1}^S \|\bX\bq_s + \tbv_s-\bone\|_2^2\\
\mathrm{over}~&~\tbz:=(\bbv,\bc,\bdelta,\bbq)\label{eq:MINLP:over}\\
\mathrm{s.to}~&~\eqref{eq:1547con2:v},\eqref{eq:1547con2:delta},\eqref{eq:1547con2:q},\eqref{eq:stability3},\eqref{eq:1547con3:cq}\label{eq:MINLP:stabilityIEEE}\\
~&~ \eqref{eq:kkt:loq}-\eqref{eq:kkt:df}\label{eq:MINLP:kkt1}&~\forall~s\\
~&~ \eqref{eq:kkt:cs1}-\eqref{eq:kkt:cs4}~\text{as in}~ \eqref{eq:bigM}&~\forall~s.\label{eq:MINLP:kkt2}
\end{align}
\end{subequations}
The bilinear terms in \eqref{eq:kkt:loq} and \eqref{eq:1547con3:cq}, and the binary variables in \eqref{eq:MINLP:kkt2} increase with the number of inverters and scenarios.

\begin{remark}\label{remark:parameterizations}
The Volt/VAR curve of Fig.~\ref{fig:curve} has four degrees of freedom that control the center, deadband, slope, and saturation of the curve. These degrees of freedom are amenable to different equivalent parameterizations, such as $(\bbv,\balpha,\bdelta,\bsigma)$ and $(\bbv,\bc,\bdelta,\bsigma)$ that we used in Sec.~\ref{sec:1pDNN}; $(\bbv,\balpha,\bdelta,\bbq)$; or $(\bbv,\bc,\bdelta,\bbq)$. We used the last one in \eqref{eq:MINLP:over} as it yielded significantly shorter solution times during our tests.
\end{remark}

% For incremental rules, the MINLP in \eqref{eq:MINLP} can be further simplified as the IEEE 1547 and stability constraints in \eqref{eq:MINLP:stabilityIEEE} are not needed anymore. 
Although the MINLP can solve ORD to near-global optimality (modulo the bilinear terms left to be handled internally by the solver), it was found to scale unfavorably with the number of DERs and/or scenarios of Section~\ref{sec:tests}. %Nonetheless, it was developed and tested to serve as a benchmark for DNN-based designs.

%%%%%%%%%%%%%%%%%%%%%%%%%%%%%%%%%%%
\section{ORD for $3\phi$ Feeders via Deep Learning}\label{sec:3pDNN}
Under transposed lines and balanced injections, one could deal with ORD using the single-phase formulations discussed earlier. Under imbalance conditions, however, a linearized multiphase feeder model would be a better approximation. DERs would still implement local Volt/VAR rules, yet sensitivity matrix $\bX$ now has different properties as discussed in \eqref{eq:ldf}.
%%%%%%%%%%%%%%%%%%%%%%%%%%%%%%%%%%%%%%%%%%%%%%%%%%%%%%%%%%
% \subsection{Equilibrium and Stability of $3\phi$}
For the multiphase case, we were not able to come up with an optimization problem whose minimizer coincides with the equilibrium $\bq^*$ similar to \eqref{eq:inner}. Nonetheless, we show in the appendix that the Volt/VAR rules of Fig.~\ref{fig:curve} do converge to a unique equilibrium under the following polytopic conditions, which form a restriction of~$\|\diag(\balpha)\bX\|_2\leq{1-\epsilon}$.

\begin{proposition}\label{pro:contraction}
Consider the Volt/VAR dynamics of \eqref{eq:dynamics1} operating over a multiphase feeder modeled by \eqref{eq:ldf}. If the Volt/VAR slope vector $\balpha$ satisfies
\begin{subequations}\label{eq:stability_3p}
\begin{align}
|\bX|^{\top}\balpha &\leq (1-\epsilon_1)\cdot \bone\label{eq:stability_3p:C1}\\
0\leq \alpha_n&\leq \frac{{1-\epsilon_2}}{\sum_{m\in\mcN}|X_{nm}|},\quad\forall n\in\mcN\label{eq:stability3p:C2}
\end{align}
for some $\epsilon_1$, $\epsilon_2$, and $\epsilon\in(0,1)$ with $(1-\epsilon_1)(1-\epsilon_2)\leq (1-\epsilon)^2$, the dynamics in~\eqref{eq:dynamics1} exhibit a unique equilibrium $\bq^*$ to which they converge exponentially fast as
\begin{align}
\|\bq^{t}-\bq^*\|_2\leq 2\|\hbq\|_2 \cdot (1-\epsilon)^t.\label{eq:conv3p}
\end{align}
\end{subequations}
\end{proposition}

The absolute value $|\bX|$ applies entry-wise. The result generalizes \eqref{eq:stability2} and \cite[Th.~3]{9091863} to multiphase feeders, wherein $\bX$ is non-symmetric and with some of its entries being negative. It provides linear constraints on $\balpha$ to ensure stability. Although $\epsilon_1$ and $\epsilon_2$ could be selected to reduce conservatism of the restriction, they were henceforth set equal as $\epsilon_1=\epsilon_2=\epsilon$.

The ORD task for multiphase feeders can be formulated as in \eqref{eq:inner} with the appropriate modification of the sensitivity matrix $\bX$. Since equilibrium setpoints cannot be expressed as the minimizer of an inner optimization, the MINLP approach of Section~\ref{sec:1pMINLP} cannot be adopted here. Alternatively, one may pursue an MINLP formulation along the lines of~\cite{STKSL22}, though scalability is still expected to be an issue. Fortunately, the DNN-based approach for ORD remains applicable with the next minor modifications: \emph{i)} Sensitivity matrices are modified accordingly; \emph{ii)} Every layer now consists of $3N$ building modules corresponding to bus/phase (node) combinations; and \emph{iii)} Use the stability constraints of \eqref{eq:stability_3p} instead of \eqref{eq:stability2}.

Proposition~\ref{pro:T} on minimum depth $T$ of DNNs for Volt/VAR rules in single-phase feeders carries over to multiphase feeders. This is easily confirmed by applying the steps from the proof of Proposition~\ref{pro:T} to the results from Proposition~\ref{pro:contraction}.

\begin{remark}\label{re:stability}
The conditions of Prop.~\ref{pro:contraction} are general enough to ensure the stability of Volt/VAR rules in any type of distribution network, single-phase or multi-phase; radial or meshed. 
\end{remark}

%%%%%%%%%%%%%%%%%%%%%%%%%%%%%%%%%%%%%%%%%%%%%%%%%%%%%%%

\section{Numerical Tests}\label{sec:tests}
The proposed ORD methods were evaluated on single- and multi-phase feeders. Real-world data of active load and solar generation at one-minute frequency was sourced from the Smart* project on April 2, 2011~\cite{Smartsolar}. The set consists of active loads from 444 homes and generation from $43$ solar panels. Loads from multiple homes were averaged to better simulate loads at buses of the primary distribution network. Each averaged load was normalized so its peak value during the day coincided with the nominal active power load of its hosting node. For each time interval, reactive loads were added by randomly sampling lagging power factors within $[0.9,1]$. Similarly to loads, each solar generation signal was normalized so its peak value was twice that of the nominal active load of the hosting bus. Apparent power limits for inverters were set to $1.1$ times the peak active generation.

Control rules were designed and evaluated in Python on a 2.4 GHz 8-Core Intel Core i9 processor laptop with 64 GB RAM. Pytorch was selected as the library to design and train DNNs, as it implements computation graphs dynamically. That is quite important for our purposes, as dynamic computation graphs imply that the number of layers $T$ does not need to be fixed beforehand. It is rather decided on the fly based on the convergence of rules for the given $\tilde{\bv}_t$. This flexibility enables limiting the DNN to lower depths. Convergence was determined based on the change in objective value by adding a layer: The rules were assumed to have converged if the objective changed by less than $1\cdot10^{-7}$ within consecutive layers. The batch size $B$ was set to the maximum of the integer part of $S/10$ and $1$. The step size $\mu$ was determined individually for each network such that the proposed design approach worked for various times of the day. All DNNs were trained using the Adam optimizer.%~\cite{Adam}

The projection step \eqref{eq:project} was performed by solving a SOCP using the CVXPY library in Python with GUROBI. The MINLP \eqref{eq:MINLP} was implemented in MATLAB using YALMIP~\cite{YALMIP} with GUROBI, and used to benchmark the results for optimality and runtime. Other details such as learning rates for DNN training, initialization of design parameters, load and solar panel assignments, and time period for scenario sampling are presented along with the corresponding results. 

%%%%%%%%%%%%%%%%%%%%%%%%%%%%%%%%%%
\subsection{Tests on Single-Phase Feeder}
\begin{figure}
	\centering
	\includegraphics[scale=0.55]{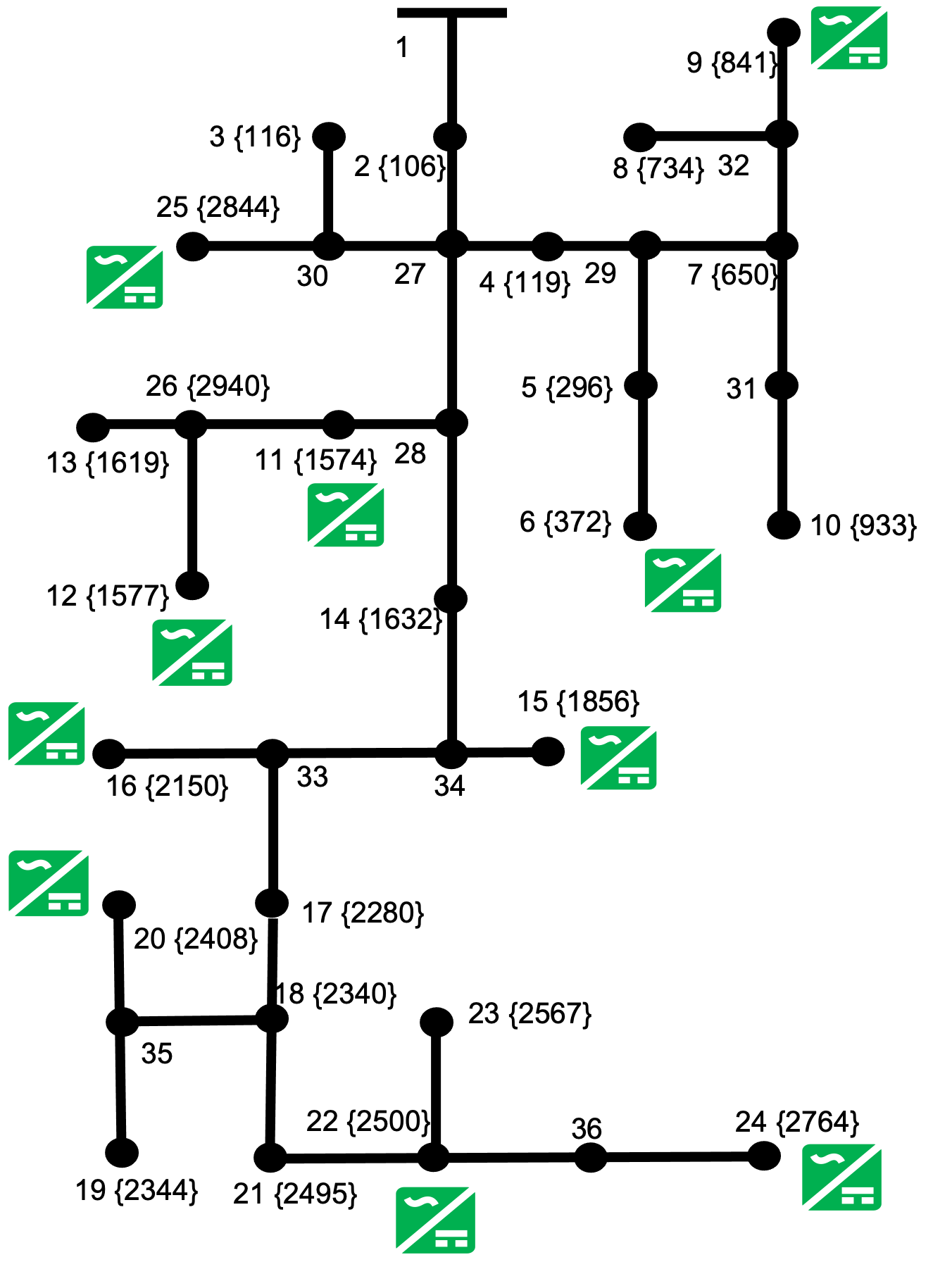}
% 	\vspace*{-1.5em}
	\caption{The IEEE 37-bus feeder used for the tests. Node numbering follows the format \texttt{node number \{panel ID\}}. DERs at nodes $\{6,9,11,12,15,16,20,22,24,25\}$ provide reactive power control; the rest operate at unit power factor.}
% 	\vspace*{-1.5em}
	\label{fig:37bus}
\end{figure}

The first set of tests was conducted on the single-phase equivalent of the IEEE 37-bus feeder. Homes with IDs $20$-$369$ were averaged $10$ at a time and successively added as active loads to buses $2-26$ as shown in Fig.~\ref{fig:37bus}. Active generation from solar panels was also added, as per the mapping in Fig.~\ref{fig:37bus}. Additionally, buses $\{6,9,11,12,15,16,20,22,24,25\}$ were equipped with DERs capable of reactive power control.

\begin{figure}[t]
	\centering
	\includegraphics[scale=0.45]{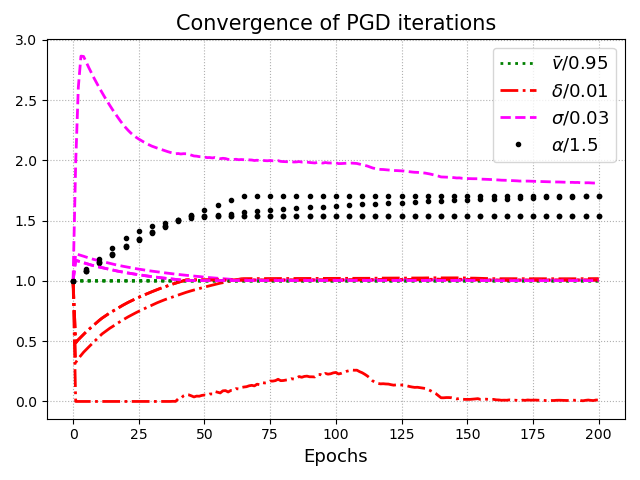}
 	% \vspace*{-1.5em}
	\caption{Convergence of PGD iterations \eqref{eq:spgd1}--\eqref{eq:project} for Volt/VAR rules with $\epsilon=0.5$. Values of rule parameters for DERs 12, 22, and 29 are plotted against training epochs. Plots have been normalized with respect to their initial values.}
 	% \vspace*{-1.5em}
	\label{fig:pgd_ni}
\end{figure}

The DNN-based rules were optimized using $80$ grid scenarios sampled from the high-solar period 3:00-4:20~pm, and were trained using the learning rate $\mu=0.003$ over $200$ epochs. The design parameters $\bz:=(\bbv,\bdelta,\bsigma,\balpha)$ were initialized at the feasible point $(\bar{v}_n,\delta_n,\sigma_n,\alpha_n)=(0.95,0.1,0.3,1.5)$ for all $n$. Figure~\ref{fig:pgd_ni} shows the convergence of Volt/VAR rule parameters for DERs at nodes $\{12,22,29\}$, for $\epsilon=0.5$, during training. To accommodate different ranges of magnitudes, all plots are normalized with respect to their initial values. 

\begin{figure}[t]
	\centering
	\includegraphics[scale=0.5]{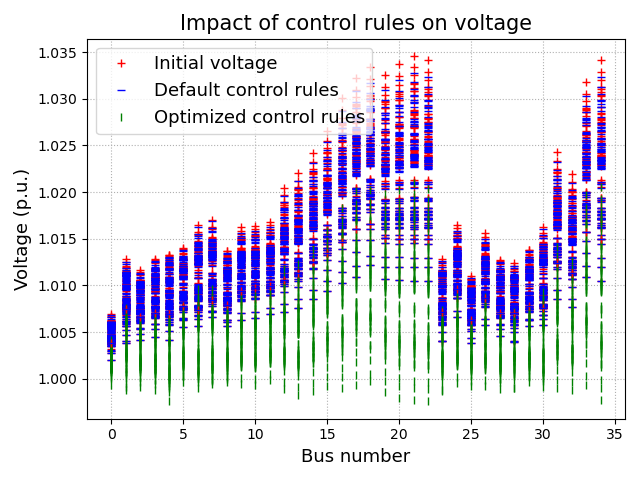}
 	% \vspace*{-1em}
	\caption{Voltage profiles across different scenarios and control options: This plot shows the voltage magnitude across all buses. Each point corresponds to a different grid loading scenario. The three colors correspond to different control options: \emph{Red} denotes no reactive power compensation by DERs (unit power factor): \emph{Blue} denotes Volt/VAR rules with the default rule parameters suggested by the IEEE Std. 1547; and \emph{Green} denotes optimized Volt/VAR rules.}
 	% \vspace*{-0.5em}
	\label{fig:profiles}
\end{figure}

Figure~\ref{fig:profiles} highlights the efficacy of the optimized Volt/VAR rules in improving the voltage profile across the feeder. Voltages across buses are plotted under three setups: voltages without DER reactive power support, voltages under the default settings $(\bar{v}_n,\delta_n, \sigma,\overline{q}_n)=(1,0.02, 0.08,\hat{q}_n)$ from IEEE 1547.8~\cite{IEEE1547}; and voltages under control rules with optimal $\bz$. For each bus, voltages for all $S=80$ scenarios have been marked. The default control rules were found only marginally to improve voltage profiles. On the other hand, optimally designed control rules successfully lowered voltages and brought them close to unity on all buses.

We next studied the impact of the stability margin $\epsilon$ on the optimal cost $L(\bz)$ of \eqref{eq:outer_DNN} under Volt/VAR rules. Recall that $\epsilon$ determines the feasible space of design parameters via \eqref{eq:stability2}. The larger the $\epsilon$, the more restricted problem~\eqref{eq:outer} is. Table~\ref{tab:epsilon} confirms this by presenting the objectives during training for a range of $\epsilon$ values. Table~\ref{tab:epsilon} also lists the chosen initial value for $\balpha$, represented by $\balpha_{\mathrm{init}}$, that renders the initial $\bz$ feasible for the corresponding value of $\epsilon$. The objective converged to the highest value for $\epsilon=0.9$ and the lowest for $\epsilon=0.5$. Note that for the studied scenarios, reducing $\epsilon$ below $0.5$ did not impact the optimal value of the objective, which indicates that the feasible space for $\epsilon=0.5$ contains the optimizers for all  $\epsilon\leq 0.5$ as well. Consequently, $\epsilon$ has been fixed at $0.5$ for the subsequent results on the $37$-bus feeder.

\begin{table}[t]
\centering
\caption{\label{tab:epsilon} Test results capturing the effect of $\epsilon$ on the optimal objective value for Volt/VAR control rules. The smaller the $\epsilon$, the larger the feasible region for rule parameters is, and so lower voltage regulation values can be attained.}
\begin{tabular}{|c|c|c|}
\hline
$\epsilon$ & $\balpha_{\mathrm{init}}$ & Objective (p.u.)\\
\hline
$0.9$ & $0.4$ & $2.22\cdot10^{-3}$\\
\hline
$0.8$ & $0.5$ & $1.37\cdot10^{-3}$\\
\hline
$0.7$ & $1$ & $1.06\cdot10^{-3}$\\
\hline
$0.6$ & $1.5$ & $9.73\cdot10^{-4}$\\
\hline
$0.5$ & $1.5$ & $8.50\cdot10^{-4}$\\
\hline
\end{tabular}
\end{table}

\begin{table}[t]
\centering
\caption{\label{tab:MINLP_S} Tests comparing the MINLP with the DNN-based ORD for different numbers of scenarios $S$, with $N_G=5$ smart DERs.}
% \vspace*{-1em}
\begin{tabular}{|c|ccc|cc|}
\hline
\multirow{2}{*}{$S$} & \multicolumn{3}{c|}{MINLP}                                                         & \multicolumn{2}{c|}{DNN}                         \\ \cline{2-6} 
                           & \multicolumn{1}{c|}{Solved} & \multicolumn{1}{c|}{Time (s)} & Obj. (p.u.)                                 & \multicolumn{1}{c|}{Time (s)} & Obj. (p.u.) \\ \hline
$10 $                        & \multicolumn{1}{c|}{Yes}    & \multicolumn{1}{c|}{$2.45$}     &\multicolumn{1}{c|}{ $9.82\cdot10^{-4}$} & \multicolumn{1}{c|}{$18.16$}      &      \multicolumn{1}{c|}{$ 9.82\cdot10^{-4}$}        \\ \hline
$20  $                       & \multicolumn{1}{c|}{Yes}    & \multicolumn{1}{c|}{$3.64$}     & \multicolumn{1}{c|}{ $1.57\cdot10^{-3}$}                                                 & \multicolumn{1}{c|}{$20.08$}         &   \multicolumn{1}{c|}{ $1.58\cdot10^{-3}$}               \\ \hline
$40  $                       & \multicolumn{1}{c|}{Yes}     & \multicolumn{1}{c|}{$123.32$}    &   \multicolumn{1}{c|}{ $2.78\cdot10^{-3}$}                                               & \multicolumn{1}{c|}{$20.63$}         &  \multicolumn{1}{c|}{ $2.78\cdot10^{-3}$}     \\ \hline
$80$                         & \multicolumn{1}{c|}{No}     & \multicolumn{1}{c|}{$300$}    &      \multicolumn{1}{c|}{ $2.68\cdot10^{-3}$}                                            & \multicolumn{1}{c|}{$22.17$}         &   \multicolumn{1}{c|}{ $2.62\cdot10^{-3}$}                  \\ \hline
\end{tabular}
% \vspace*{-1em}
\end{table}

\begin{figure}[t]
	\centering
	\includegraphics[scale=0.5]{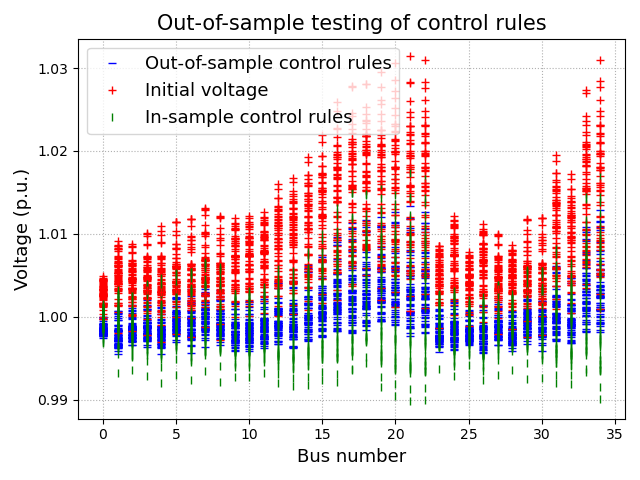}
 	% \vspace*{-1em}
	\caption{Bus voltages of the IEEE 37-feeder across $S=60$ scenarios during 2:00-3:00~pm. \emph{In-sample} rules were designed and validated using the same $S=60$ scenarios drawn between 2:00-3:00~pm. \emph{Out-of-sample} rules were designed using $S=60$ scenarios from 1:00-2:00~pm, and validated using the $S=60$ scenarios drawn from 2:00-3:00~pm. Voltages experienced during 2:00--3:00~pm with no Volt/VAR control are also shown for comparison. As expected, the performance of the rules degrades when rules are designed based on non-representative scenarios.}
 	% \vspace*{-0.5em}
	\label{fig:oos_ac}
\end{figure}

We also studied the performance of control rules on scenarios other than those used for designing the rules. To obtain different loading conditions over time, we scaled the overall solar generation by $0.8$. Then two sets of control rules were designed. The first set of rules was designed using $S=60$ samples drawn between 2:00-3:00~pm and is referred to as the \emph{in-sample} sample rules. The second set of rules was designed using $S=60$ samples drawn between 1:00-2:00~pm and is referred to as the \emph{out-of-sample} control rules. Both sets of rules were then evaluated on the samples drawn between 2:00-3:00~pm. Figure~\ref{fig:oos_ac} illustrates the voltage experienced on each bus and across all $S=60$ scenarios for the two sets of control rules. The voltages experienced with no Volt/VAR control are also plotted for comparison. Both sets of rules improved the voltage profile over the no Volt/VAR control option, yet out-of-sample exhibited a somewhat inferior performance compared to the in-sample rules as expected.

\begin{figure}[t]
	\centering
	\includegraphics[scale=0.5]{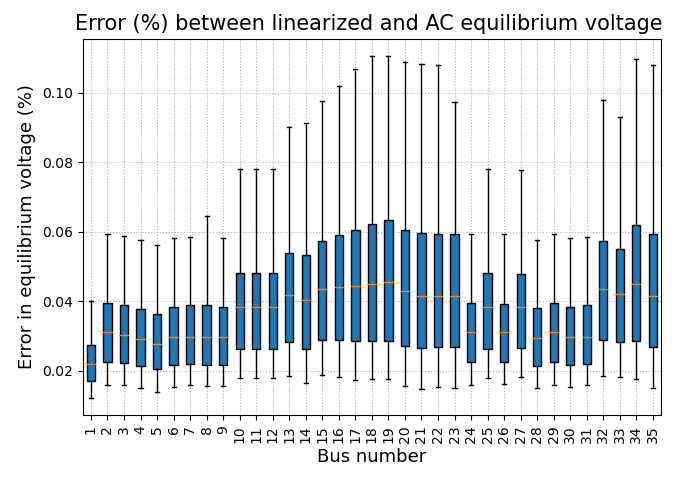}
 	% \vspace*{-1em}
	\caption{Comparison of equilibrium voltages under the linearized and the exact AC grid models. Control rules were designed over the 2:00--3:00~pm  window using the linearized model and then ran until reaching equilibrium using the linearized and AC models. Errors between the equilibrium voltages, captured in the percentage of AC equilibrium voltage, were found to be small enough to validate the proposed approach for designing control rules.}
 	% \vspace*{-0.5em}
	\label{fig:ac_error}
\end{figure}

While rules were designed using the linearized grid model, their performance on the accurate AC grid model was also evaluated. For this purpose, the \emph{in-sample} control rules designed in the previous paragraph for 2:00--3:00~pm  were applied and ran until equilibrium, under the linearized and AC power flow models. Then for each bus and scenario, a percentage error was calculated by taking the difference between the linearized and corresponding AC equilibrium voltage, divided by the AC equilibrium voltage. The results are illustrated as a box plot in Fig.~\ref{fig:ac_error}. Evidently, the linearization error at equilibrium is consistently less than 0.1\% or approximately 0.001 per unit, which verifies that control rules are effective over the AC grid model too.

 \begin{figure}[t]
	\centering
	\includegraphics[scale=0.5]{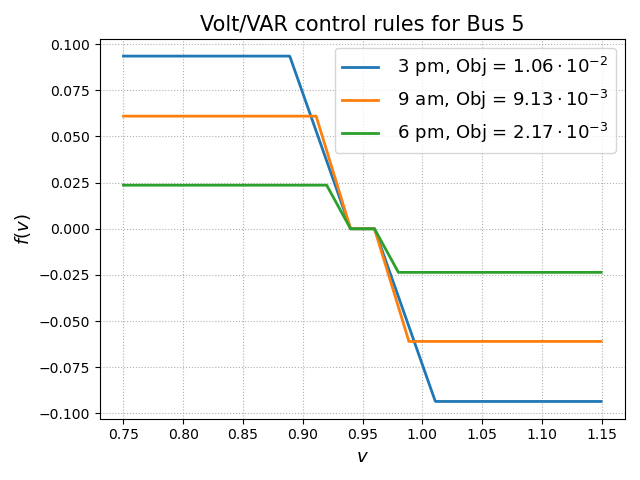}
 	% \vspace*{-1em}
	\caption{Volt/VAR rules obtained for Bus 5. Three sets of voltage scenarios were sampled. The deviation from the nominal 1~pu is captured by the value for \texttt{Obj} for each of the sets. The Volt/VAR rules become steeper with higher saturation limits as the grid voltages worsen.}
 	% \vspace*{-0.5em}
	\label{fig:control_curves_many}
\end{figure}

The last test on the IEEE 37-bus feeder intended to provide some intuition on the shape of the control rules designed across different periods of the day. To this end, we sampled three sets of $S=80$ scenarios at three different hours of the day: 9~am, 3~pm, and 6~pm. For each set, we quantified the experienced voltage profiles by a single number given by the objective of~\eqref{eq:outer}, before voltage regulation. We fixed the design parameters $\bbv$ and $\bdelta$ to $0.95$ and $0.01$, respectively, hence localizing the impact of the control rule design process to $\bsigma$ and $\balpha$. The resulting Volt/VAR rules for Bus~5 are presented in Figure~\ref{fig:control_curves_many}. The legend for each rule captures the starting time for drawing $S=80$ consecutive samples. The value of the objective of~\eqref{eq:outer} for prior voltages is also presented for each rule. Figure~\ref{fig:control_curves_many} shows that as the grid voltages move further away from 1~pu, both $\balpha$ and $\bsigma$ increase, resulting in steeper curves with higher saturation limits. Such a design enables the curves to provide more voltage regulation as grid voltages worsen.

To verify the optimality and scalability of DNN-based ORD, we benchmarked them against the MINLP formulation of~\eqref{eq:MINLP}. The MINLP was allowed to run until completion or till $300$ seconds, whichever happened earlier. Scaling with respect to both the number of scenarios as well as DERs was studied. Table~\ref{tab:MINLP_S} reports the results for the case when the number of smart DERs was fixed to $N_G=5$ and scenarios were increased from $S=10$ to $80$. As evident from Table~\ref{tab:MINLP_S}, the DNN-based ORD scaled much better than the MINLP for larger $S$, as expected. Furthermore, the DNN-based ORD achieved the same objective as the MINLP across all tested values of $S$. This is remarkable since SGD for non-convex problems can only guarantee convergence to stationary points. Similar conclusions can be drawn from Table~\ref{tab:MINLP_I} where we fixed $S=80$ and varied $N_G$ from $2$ to $10$. The MINLP was faster than the DNN-based approach for $N_G = 2$, but could not be solved within $300$ seconds if more inverters were added. On the other hand, the DNN-based ORD scaled gracefully with the $N_G$ and achieved lower objectives for all $N_G\geq4$.

\begin{table}[t]
\centering
\caption{\label{tab:MINLP_I} Comparing MINLP with the DNN-based approach\\ for different $N_G$ and $S=80$.}
% \vspace*{-1em}
\begin{tabular}{|c|ccc|cc|}
\hline
\multirow{2}{*}{$N_G$} & \multicolumn{3}{c|}{MINLP}                                                         & \multicolumn{2}{c|}{DNN}                         \\ \cline{2-6} 
                           & \multicolumn{1}{c|}{Solved} & \multicolumn{1}{c|}{Time (s)} & Obj. (p.u.)                                 & \multicolumn{1}{c|}{Time (s)} & Obj. (p.u.) \\ \hline
$2 $                        & \multicolumn{1}{c|}{Yes}    & \multicolumn{1}{c|}{$3.90$}     &\multicolumn{1}{c|}{ $3.62\cdot10^{-3}$} & \multicolumn{1}{c|}{$14.12$}      &      \multicolumn{1}{c|}{$ 3.62\cdot10^{-3}$}        \\ \hline
$4  $                       & \multicolumn{1}{c|}{No}    & \multicolumn{1}{c|}{$300$}     & \multicolumn{1}{c|}{ $3.22\cdot10^{-3}$}                                                 & \multicolumn{1}{c|}{$17.96$}         &   \multicolumn{1}{c|}{ $3.18\cdot10^{-3}$}               \\ \hline
$6  $                       & \multicolumn{1}{c|}{No}     & \multicolumn{1}{c|}{$300$}    &   \multicolumn{1}{c|}{ $2.77\cdot10^{-3}$}                                               & \multicolumn{1}{c|}{$21.95$}         &  \multicolumn{1}{c|}{ $2.35\cdot10^{-3}$}     \\ \hline
$8$                         & \multicolumn{1}{c|}{No}     & \multicolumn{1}{c|}{$300$}    &      \multicolumn{1}{c|}{ $1.40\cdot10^{-3}$}                                            & \multicolumn{1}{c|}{$33.42$}         &   \multicolumn{1}{c|}{ $1.16\cdot10^{-3}$}                  \\ \hline
$10$                         & \multicolumn{1}{c|}{No}     & \multicolumn{1}{c|}{$300$}    &      \multicolumn{1}{c|}{ $1.20\cdot10^{-3}$}                                            & \multicolumn{1}{c|}{$39.76$}         &   \multicolumn{1}{c|}{ $8.50\cdot10^{-4}$}                  \\ \hline
\end{tabular}
\end{table}

\begin{figure}[t]
	\centering
	\includegraphics[scale=0.4]{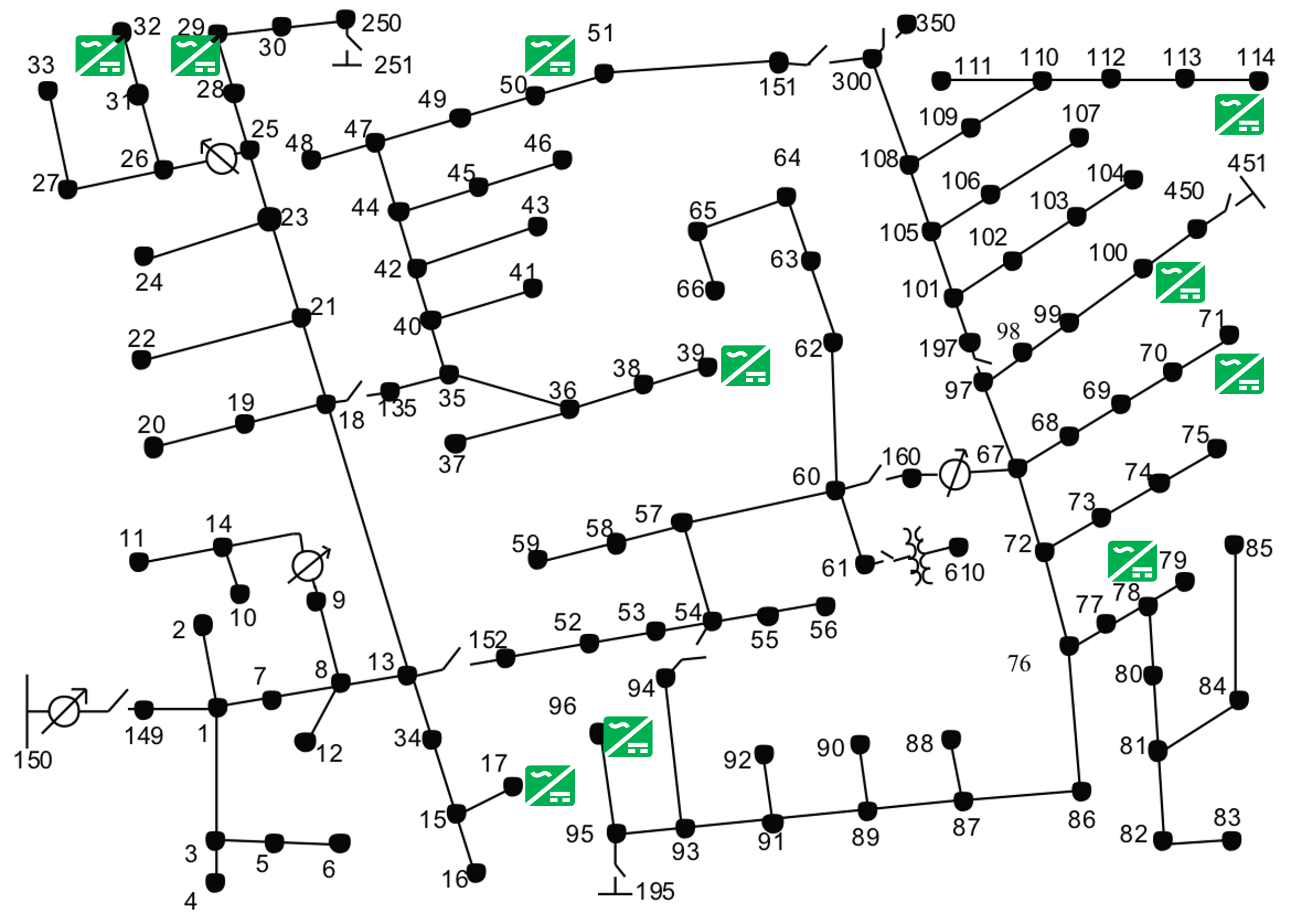}
% 	\vspace*{-1.5em}
	\caption{Inverter siting on the IEEE 123-bus distribution feeder.}
% 	\vspace*{-1.5em}
	\label{fig:ieee123}
\end{figure}

The scalability of the DNN-based control rules was also confirmed by implementing them for the larger IEEE 123-bus feeder of Fig.~\ref{fig:ieee123}. Active load data was generated by averaging homes with IDs $20$-$386$, three at a time, and were serially assigned to buses $2$-$123$. Solar generation from $10$ panels with IDs $\{106,116,119,296,372,650,734,841,933,1574\}$ was added to buses $\{17,29,32,39,50,71,78,96,100,114\}$, respectively. All solar buses were equipped with smart DERs for reactive power support. The DNNs for Volt/VAR rules were trained with the learning rate  $\mu=0.01$, with $\epsilon$ set to $0.5$. The design parameters $\bz:=(\bbv,\bdelta,\bsigma,\balpha)$ were initialized at the feasible point $(1.05,0.1,0.3,1.5)$, and $\mu$ was set as $0.01$. With $N_G$ and $S$ fixed at $10$ and $80$, respectively, the DNN-based ORD was compared to the MINLP one. For this larger network, the MINLP solver was allowed to run until $500$ seconds. To ensure repeatability, the results were repeated across several time periods between $1-6$~PM, and have been compiled in Table~\ref{tab:MINLP_123}. For all time periods, the DNN-based solver scaled well in terms of the DNN training time. The MINLP solver could not converge within $500$ seconds and was outperformed by the DNN-based solver in terms of final objective values across all setups.

\begin{table}[t]
\centering
\caption{\label{tab:MINLP_123} Test results comparing the MINLP with the DNN-based ORD approach for the single-phase IEEE 123-bus feeder, across different time periods for $N_G=10$ DERs and $S=80$ scenarios.}
% \vspace*{-1em}
\begin{tabular}{|c|ccc|cc|}
\hline
\multirow{2}{*}{Time} & \multicolumn{3}{c|}{MINLP}                                                         & \multicolumn{2}{c|}{DNN}                         \\ \cline{2-6} 
                           & \multicolumn{1}{c|}{Solved} & \multicolumn{1}{c|}{Time (s)} & Obj. (p.u.)                                 & \multicolumn{1}{c|}{Time (s)} & Obj. (p.u.) \\ \hline
$1$ pm                        & \multicolumn{1}{c|}{No}    & \multicolumn{1}{c|}{$500$}     &\multicolumn{1}{c|}{ $9.26\cdot10^{-4}$} & \multicolumn{1}{c|}{$28.6$}      &      \multicolumn{1}{c|}{$ 8.95\cdot10^{-4}$}        \\ \hline
$2$ pm                       & \multicolumn{1}{c|}{No}    & \multicolumn{1}{c|}{$500$}     & \multicolumn{1}{c|}{ $6.69\cdot10^{-4}$}                                                 & \multicolumn{1}{c|}{$30.18$}         &   \multicolumn{1}{c|}{ $6.40\cdot10^{-4}$}               \\ \hline
$3$ pm                       & \multicolumn{1}{c|}{No}     & \multicolumn{1}{c|}{$500$}    &   \multicolumn{1}{c|}{ $4.17\cdot10^{-4}$}                                               & \multicolumn{1}{c|}{$27.55$}         &  \multicolumn{1}{c|}{ $3.92\cdot10^{-4}$}     \\ \hline
$4$ pm                        & \multicolumn{1}{c|}{No}     & \multicolumn{1}{c|}{$500$}    &      \multicolumn{1}{c|}{ $2.17\cdot10^{-4}$}                                            & \multicolumn{1}{c|}{$29.83$}         &   \multicolumn{1}{c|}{ $2.09\cdot10^{-4}$}                  \\ \hline
$5$ pm                        & \multicolumn{1}{c|}{No}     & \multicolumn{1}{c|}{$500$}    &      \multicolumn{1}{c|}{ $2.98\cdot10^{-3}$}                                            & \multicolumn{1}{c|}{$27.53$}         &   \multicolumn{1}{c|}{ $2.87\cdot10^{-4}$}                  \\ \hline
\end{tabular}
\end{table}

%%%%%%%%%%%
\subsection{Tests on a Multiphase Feeder}
\begin{figure}[t]
	\centering
	\includegraphics[scale=0.45]{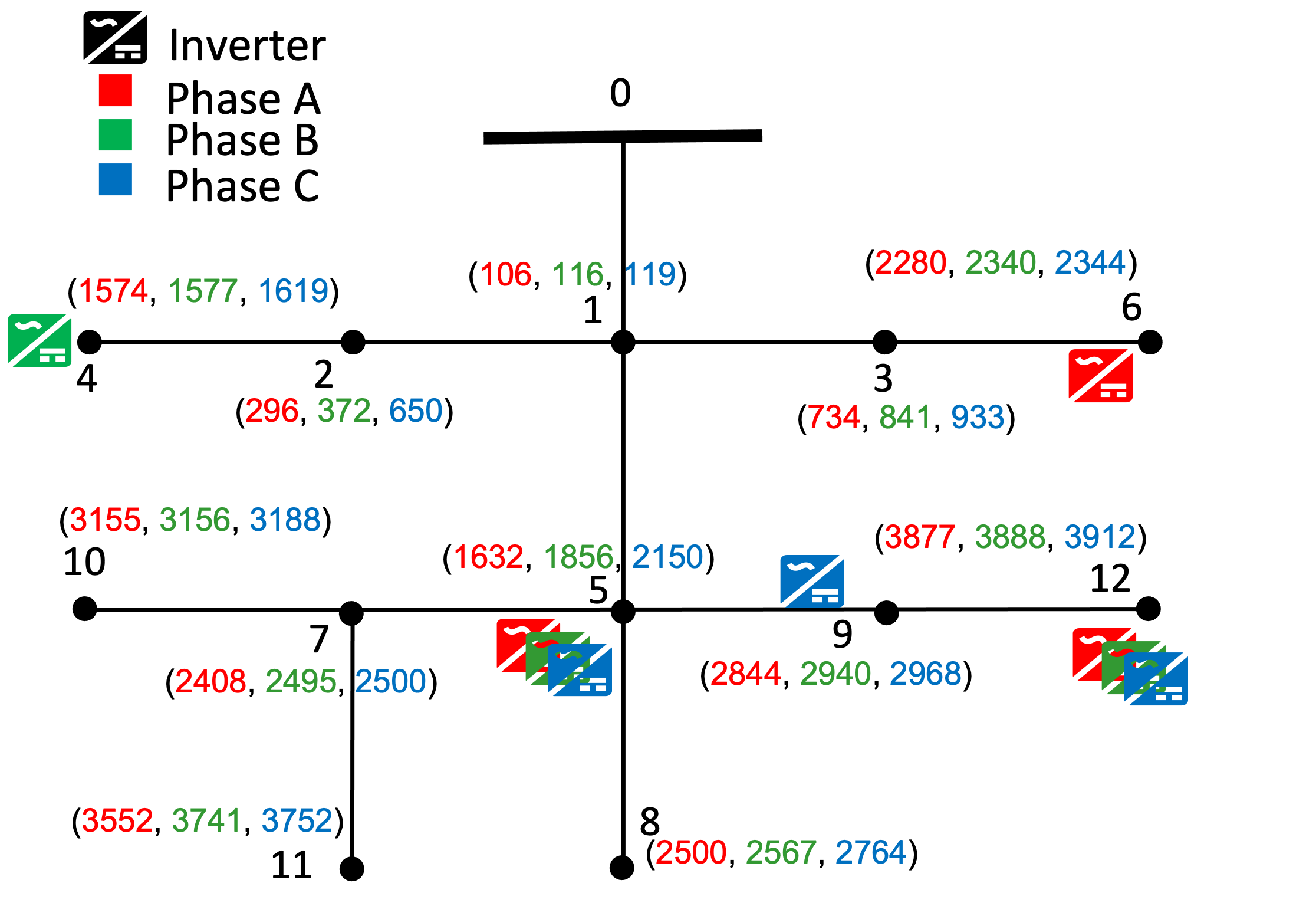}
% 	\vspace*{-1.5em}
	\caption{The three-phase IEEE 13-bus distribution feeder system.}
% 	\vspace*{-1.5em}
	\label{fig:ieee13}
\end{figure}

The DNN-based control rules were also tested on the multiphase IEEE 13-bus feeder. Active loads from homes with IDs $20$-$379$ were averaged ten homes at a time. The resulting $36$ averaged loads were added to buses $1$-$12$, allocating all three phases for a bus before moving on to the next one. Solar generation was added to nodes per the panel assignments in Fig.~\ref{fig:ieee13}. Values in red, green, and blue correspond to panel IDs assigned to Phases A, B, and C, respectively. Reactive power compensation was provided by nine inverters added across phases, and bus indices, as shown in Fig~\ref{fig:ieee13}, with the colors indicating the corresponding phase. 

The learning rate for DNN-based control rules was set to $\mu=0.1$, with the design parameters $\bz:=(\bbv,\bdelta,\bsigma,\balpha)$ initialized to feasible values $(0.95,0.01,0.3,1.5)$. In the absence of an MINLP solver, the optimized DNN-based control rules were benchmarked against control rules with the default settings from the IEEE 1547.8 standard. Table~\ref{tab:13_nivsi} collects the values for the objective~\eqref{eq:outer} for $S=80$ scenarios, across different windows of time from $1-5$ pm, under three control schemes-- no reactive power compensation, optimized control rules, and default control rules. The default control rules did not manage to significantly reduce the objective~\eqref{eq:outer}, as the grid conditions $\tbv$ were observed to fall in the deadband of the default control rules frequently. In contrast, the optimized control rules took the grid conditions $\tbv$ into consideration while designing the deadband, and hence improved voltage profiles considerably.

\begin{table}[t]
\centering
\caption{\label{tab:13_nivsi} Test results on the multiphase IEEE 13-bus feeder for $N_G=9$ inverters and $S=80$ scenarios. Comparing the objective~\eqref{eq:outer} under three scenarios: no reactive power compensation, optimized control rules, and the default rules per IEEE 1547.}
% \vspace*{-1em}
\begin{tabular}{|c|c|c|c|}
\hline
\multirow{1}{*}{Time} & \multicolumn{1}{c|}{$\bq=0$} & \multicolumn{1}{c|}{Optimized}      & \multicolumn{1}{c|}{Default}                       \\ \hline

 $1$ pm             & \multicolumn{1}{c|}{$2.51\cdot10^{-3}$} &  \multicolumn{1}{c|}{$1.15\cdot10^{-3}$}        &  $2.31\cdot10^{-3}$ \\ \hline
 $2$ pm             & \multicolumn{1}{c|}{$1.48\cdot10^{-3}$}  & \multicolumn{1}{c|}{$6.89\cdot10^{-4}$}         & $1.42\cdot10^{-4}$ \\ \hline
    $3$ pm             & \multicolumn{1}{c|}{$6.89\cdot10^{-4}$} &  \multicolumn{1}{c|}{$4.94\cdot10^{-4}$}        & $6.89\cdot10^{-4}$ \\ \hline
    $4$ pm             & \multicolumn{1}{c|}{$8.03\cdot10^{-4}$} & \multicolumn{1}{c|}{$5.26\cdot10^{-4}$}        & $8.03\cdot10^{-4}$ \\ \hline
    $5$ pm             & \multicolumn{1}{c|}{$5.51\cdot10^{-4}$}  & \multicolumn{1}{c|}{$4.11\cdot10^{-4}$}        & $5.51\cdot10^{-4}$ \\ \hline
\end{tabular}
% \vspace*{-1em}
\end{table}

% \begin{figure}[t]
% 	\centering
% % 	\includegraphics[scale=0.55]{curve.png}% 3 plots with vr fixed to 1
% 	\includegraphics[scale=0.52]{convergence_control_rules.png}
% % 	\vspace*{-1.5em}
% 	\caption{Inverter voltages under different optimized control rules for given $\tbv$ and $\bz$, plotted against the DNN layers. The non-incremental rules require the least number of DNN layers, followed by accelerated incremental control rules, and incremental control rules.}
% % 	\vspace*{-1.5em}
% 	\label{fig:conv_cr}
% \end{figure}

%%%%%%%%%%%%%%%%%%%%%%%%%%%%%%%%%%%%%%%%
\section{Conclusions}\label{sec:conclusions}
This work has genuinely reformulated the ORD problem to train a DNN using grid scenarios as training data, unit voltages as desired targets for equilibrium voltages, and Volt/VAR rule parameters as weights. The proposed DNN-based ORD framework is general enough to accommodate Volt/VAR rules on single- and multi-phase feeders. We have also reviewed and extended results on the stability and convergence rates of Volt/VAR control rules. For benchmarking purposes, we have also developed a MINLP approach to ORD. The suggested approaches have been validated using real-world data on IEEE test feeders. The tests show that DNN-based ORD outperforms the MINLP approach in terms of optimality under time budgets and that optimized ORD curves outperform the default values. Our findings form the foundations for exciting research directions, such as: \emph{d1)} Can the DNN-based ORD framework be extended to designing \emph{incremental} Volt/VAR control rules with favorable stability characteristics? \emph{d2)} What are the appropriate Volt/VAR control rules for three-phase (probably large-scale utility-owned) DERs?

\appendix
%%%%%%%%%%%% Proposition 1 %%%%%%%%%%%%%
% \begin{IEEEproof}\emph{Proposition~\ref{pro:T}:}
\emph{Proof of Proposition~\ref{pro:T}:}
By a contraction mapping argument, reference~\cite{9091863} proves that as long as stable, the Volt/VAR dynamics $\bq^t$ enjoy exponential convergence to the equilibrium $\bq^*$. That means that if $\|\diag(\balpha)\bX\|_2<1$, then
\begin{equation*}%\label{eq:conv_q}
\|\bq^{t}-\bq^*\|_2\leq \|\diag(\balpha)\bX\|_2\cdot\|\bq^{t-1}-\bq^*\|_2.  
\end{equation*}
Propagating the previous claim across time and for $\epsilon$-stable rules $\|\diag(\balpha)\bX\|_2\leq 1-\epsilon$, we get that
\begin{equation*}%\label{eq:conv_q}
\|\bq^{t}-\bq^*\|_2\leq \|\bq^0-\bq^*\|_2\cdot (1-\epsilon)^t\leq 2\|\hbq\|_2\cdot (1-\epsilon)^t
\end{equation*}
since the initial distance to the equilibrium can be upper bounded by $\|\bq^0-\bq^*\|_2\leq 2\|\hbq\|_2$. Because $\bv=\bX\bq+\tbv$, translate injection distances to voltage distances
\begin{equation*}%\label{eq:conv_q}
\|\bv^{t}-\bv^*\|_2\leq 2\|\bX\|_2\|\hbq\|_2(1-\epsilon)^t.
\end{equation*}
To ensure the voltage approximation error at time $T$ is smaller than $\epsilon_1$, or
% \begin{equation*}
$\|\bv^T-\bv^*\|_2\leq 2\|\bX\|_2\|\hbq\|_2(1-\epsilon)^T\leq \epsilon_1,$
% \end{equation*}
it suffices to select $T$ as
\begin{equation*}%\label{eq:conv_q}
T \log(1-\epsilon)\leq \log \frac{\epsilon_1}{2\|\bX\|_2\|\hbq\|_2}.
\end{equation*}
The claim follows by noticing that $\log(1-\epsilon)<0$. 
\qed
% \end{IEEEproof}

%%%%%%%%%%%% Proposition 3 %%%%%%%%%%%%%
\emph{Proof of Proposition~\ref{pro:contraction}:}
Reference \cite[Th.~3]{9091863} shows that the Volt/VAR rules of $\bef(\cdot)$ are Lipschitz continuous in $\bq$ with $\|\diag(\balpha)\bX\|_2$ as the Lipschitz constant, that is
 \begin{align}\label{eq:lipschitz}
    \|\bef(\bq)-\bef(\bq')\|_2&\leq\|\diag(\balpha)\bX\|_2\cdot \|\bq-\bq'\|_2
\end{align}
for any $\bq$ and $\bq'$ obeying~\eqref{eq:1547con2:q}. From H\"{o}lder's inequality for matrix norms, it holds that 
% \begin{subequations}%\label{eq:holders_3p}
 \begin{align*}
\|\diag(\balpha)\bX\|_2^2 &\leq \|\diag(\balpha)\bX\|_1\cdot \|\diag(\balpha)\bX\|_\infty\\
&=\|\diag(\balpha)|\bX|\|_1\cdot \|\diag(\balpha)|\bX|\|_\infty
\end{align*}
% \end{subequations}
where $\|\cdot\|_1$ and $\|\cdot\|_\infty$ are defined as the maximum absolute sums column-wise and row-wise, respectively. The equality holds because $\balpha$ has positive entries. It is easy to check that $\|\diag(\balpha)|\bX|\|_1$ is the maximum entry of vector $|\bX|^{\top}\balpha$, and $\|\diag(\balpha)|\bX|\|_\infty$ is the maximum entry of vector $\diag(|\bX|\bone)\balpha$. Consequently, enforcing \eqref{eq:stability_3p} results in $\|\diag(\balpha)\bX\|_2\leq(1-\epsilon)$. Substituting $\|\diag(\balpha)\bX\|_2<(1-\epsilon)$ in \eqref{eq:lipschitz} yields
 \begin{equation}\label{eq:contraction}
    \|\bef(\bq)-\bef(\bq')\|_2\leq(1-\epsilon)\|\bq-\bq'\|_2
\end{equation} 
Since $\epsilon\in(0,1)$, the above relation is a contraction mapping over the space $\bq\in[-\bbq,\bbq]$ with respect to the $\ell_2$-norm. The latter establishes the existence and uniqueness of the equilibrium, as well as the exponential convergence of Volt/VAR dynamics. To explicitly derive the convergence result~\eqref{eq:conv3p}, note that $\bq^{t}=\bef\left(\bq^{t-1}\right)$ and $\bq^*=\bef\left(\bq^*\right)$. From~\eqref{eq:contraction}, we get
\[\|\bq^{t}-\bq^*\|_2\leq   (1-\epsilon) \|\bq^{t-1}-\bq^*\|_2\\
  \leq(1-\epsilon)^t \|\bq^{0}-\bq^*\|_2.\]
The claim follows as $\|\bq^{0}-\bq^*\|_2\leq 2\|\hbq\|_2$.
\qed
\balance
\bibliographystyle{IEEEtran}
\bibliography{myabrv,kekatos,inverters}

\begin{IEEEbiography}[{\includegraphics[width=1in,height=1.25in,clip,keepaspectratio]{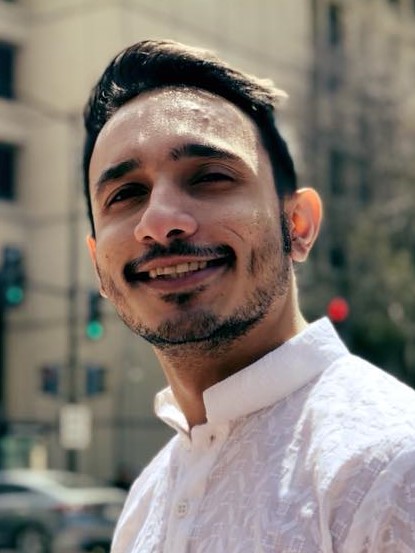}}] 
{Sarthak Gupta} received the B.Tech. degree in Electronics and Electrical Engineering from the Indian Institute of Technology Guwahati, India, in 2013, and the M.Sc. and Ph.D. degrees in Electrical Engineering from Virginia Tech, Blacksburg, VA, USA, in 2017 and 2022, respectively. In 2017 and 2019, he worked as an Associate Engineer with the Distributed Resources Operations Team, New York Independent System Operator, NY, USA, respectively. During the Summer of 2021, he interned with the Applied Mathematics and Plasma Physics Group, Los Alamos National Laboratory, Los Alamos, NM, USA. He is currently a Senior Data Scientist with C3.AI. His research interests include deep learning, reinforcement learning, optimization, and power systems.
\end{IEEEbiography}

\begin{IEEEbiography}[{\includegraphics[width=1in,height=1.25in,clip,keepaspectratio]{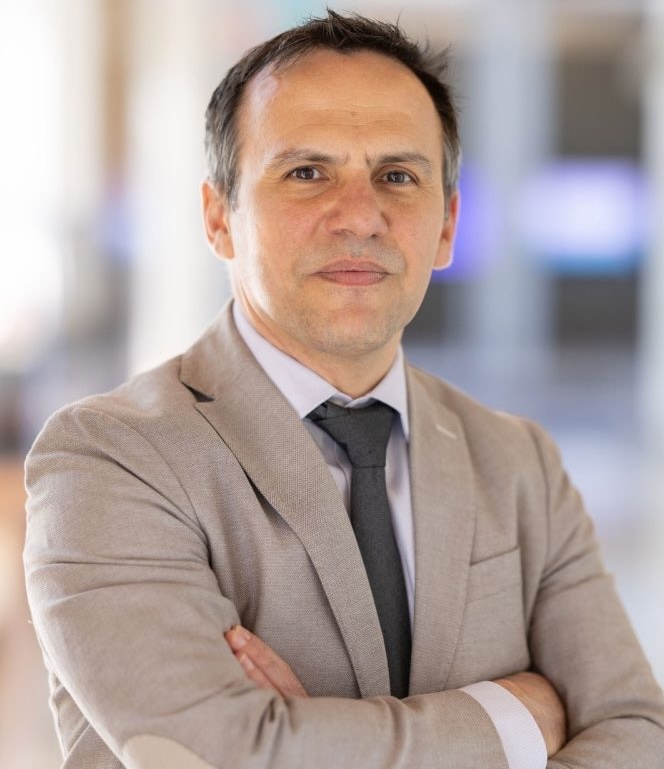}}] {Vassilis Kekatos} (SM'16) is an Associate Professor with the Schweitzer Power and Energy Systems group at the Elmore Family School of Electrical and Computer Engineering of Purdue University. He obtained his Ph.D. in Computer Science and Engineering from the Univ. of Patras, Greece in 2007. He received a Marie Curie Fellowship from the European Commission during 2009-2012, and the US National Science Foundation CAREER Award in 2018. He was a postdoctoral research associate with the ECE Dept. at the Univ. of Minnesota, and a visiting researcher with the Univ. of Texas at Austin and the Ohio State University. During 2015-2023, he was with the Bradley Department of ECE at Virginia Tech. From 2015 to 2022, he served as an Associate Editor in the IEEE Trans. on Smart Grid, and now as an Associate Editor in the IEEE Trans. on Power Systems. His current research focuses on optimization, machine learning, and quantum computing solutions for addressing power systems computational tasks.
\end{IEEEbiography}

\begin{IEEEbiography}[{\includegraphics[width=1in,height=1.25in,clip,keepaspectratio]{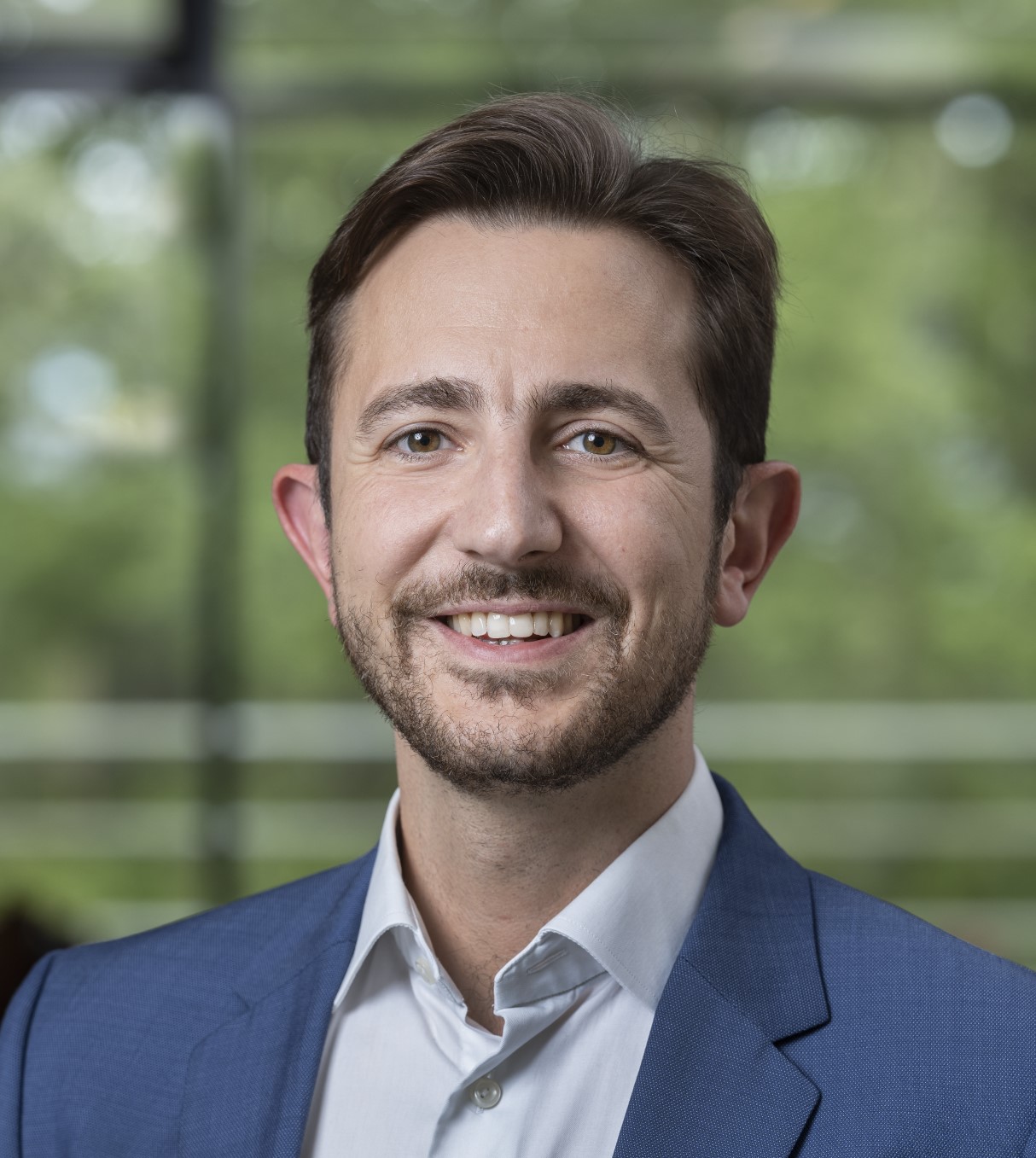}}] 
{Spyros Chatzivasileiadis} (S’04, M’14, SM’18) is a Professor and the Head of the Section for Power Systems at the Technical University of Denmark (DTU). Before that, he was a postdoctoral researcher at the Massachusetts Institute of Technology (MIT), USA, and at Lawrence Berkeley National Laboratory, USA. Spyros holds a PhD from ETH Zurich, Switzerland (2013) and a Diploma in Electrical and Computer Engineering from the National Technical University of Athens (NTUA), Greece (2007). He is currently working on trustworthy machine learning for power systems, quantum computing, and optimization, dynamics, and control of power systems. Spyros has received the Best Teacher of the Semester Award at DTU Electrical Engineering, and is the recipient of an ERC Starting Grant in 2020.
\end{IEEEbiography}

\end{document}